\documentclass[11pt,twoside,reqno]{amsproc}
\usepackage[centertags]{amsmath}
\usepackage{amssymb,amsfonts,amsthm}
\newcommand{\TryPackage}[3]{\IfFileExists{#1.sty}{\usepackage{#1} #2}{#3}}
\TryPackage{mathrsfs}{\renewcommand{\mathcal}{\mathscr}}{%
   \TryPackage{eucal}{}{}}
\usepackage{a4ML,math40,local}
\renewcommand{\subsectionmark}[1]{\relax}

\begin{document}

\title[Inverse problem for first order systems]{
The inverse spectral problem for first order systems on the half line}

\author{Matthias Lesch}
\address{Humboldt Universit\"at zu Berlin\\
Institut f\"ur Mathematik\\ 
Unter den Linden 6\\
D--10099 Berlin\\
Germany}
\curraddr{Universit\"at Bonn\\ Mathematisches Institut\\ Beringstr. 1\\
D--53115 Bonn}
\email{lesch@mathematik.hu-berlin.de, lesch@math.uni-bonn.de}
\urladdr{http://spectrum.mathematik.hu-berlin.de/$\sim$lesch}
\author{Mark Malamud}
\address{Department of Mathematics\\ University of Donetsk\\ Donetsk\\ 
Ukraine}
\email{mmm@univ.donetsk.ua}
\dedicatory{
Dedicated to the memory of M. G. Krein on the occasion of the
90th anniversary of his birth}

\begin{abstract}

On the half line $[0,\infty)$ we study
first order differential operators of the form
   \[ B\frac 1i \frac{d}{dx} + Q(x),\]
where $B:=\mat{B_1}{0}{0}{-B_2},$
$B_1,B_2\in {\rm M}(n,\C)$ are self--adjoint
positive definite matrices and $Q:\R_+\to {\rm M}(2n,\C),$
$\R_+:=[0,\infty),$ is
a continuous self--adjoint off--diagonal matrix function.

We determine the self--adjoint boundary conditions for
these operators. We prove that for each such boundary
value problem there exists a unique matrix spectral
function $\sigma$ and a generalized Fourier transform
which diagonalizes the corresponding operator 
in $L^2_{\sigma }(\R,\C)$.

We give necessary and sufficient conditions for a
matrix function $\sigma$ to be the spectral measure
of a matrix potential $Q$. Moreover we
present a procedure based on a Gelfand-Levitan type
equation for the determination of $Q$ from $\sigma $.
Our results generalize earlier results of M. Gasymov and
B. Levitan.

We apply our results to show the existence of $2n\times 2n$ Dirac
systems with purely absolute continuous, purely singular continuous
and purely discrete spectrum of multiplicity $p$, where $1\le p\le n$
is arbitrary. 
%
%
\end{abstract}

\maketitle
\tableofcontents

\section{Introduction}\label{sec0}

We consider the differential operator
\begin{equation}
    L:=B\frac 1i \frac{d}{dx} + Q(x),
   \label{G1-1.1}
\end{equation}
where 
$$B:=\mat{B_1}{0}{0}{-B_2},$$
$B_1,B_2\in {\rm M}(n,\C)$ are self--adjoint
positive definite matrices and $Q:\R_+\to {\rm M}(2n,\C),$
$\R_+:=[0,\infty),$ is
a continuous self--adjoint matrix function.
If $B_1=B_2=I_n$ then \myref{G1-1.1} is a Dirac operator.

It turns out that the operator \myref{G1-1.1} subject to the boundary condition
\begin{equation}
  f_2(0)=Hf_1(0) \quad\mbox{with}\quad B_1=H^*B_2H
  \label{G10-0.2}
\end{equation}
generates a self-adjoint extension $L_H$ of the minimal operator corresponding
to $L$. Here, $f_1(0), f_2(0)$ denote the first resp. last $n$ components
of the vector $f(0)$.

Let 
$Y(x,\lambda )$ be the $2n\times n$ matrix solution of the
initial value problem
\begin{equation}
    LY=\lambda Y,\quad  Y(0,\lambda )={I\choose H}.
    \label{G10-0.3}
\end{equation}

We will prove that there exists a unique increasing right--continuous
$n\times n$ matrix
function $\sigma (\lambda ),\lambda \in \R$, (spectral function or spectral
measure) such that we have the
symbolic identity

\begin{equation}
    \int_\R Y(x,\lambda )d\sigma (\lambda )Y(t,\lambda )^*=
\delta (x-t)I_{2n}.
    \label{G10-0.4}
\end{equation}

The main purpose of this paper is to investigate the inverse spectral problem
for the operator $L_H$. This means to find necessary and sufficient
conditions for a $n\times n$ matrix function $\sigma $ to be the spectral
function of the boundary value problem \myref{G1-1.1}, \myref{G10-0.2}.

For a Sturm-Liouville operator this problem has been posed and completely
solved by I. Gelfand and B. Levitan in the well--known paper
\cite{GelLev:DDESF}
(see also \cite{Kre:TFODSOBVP}, \cite{LevSar:SLDO},
\cite{Mar:SLOA}). 
Later on M. Gasymov and B. Levitan proved similar results
for $2\times 2$ Dirac systems \cite{GasLev:IPDS}, \cite[Chap. 12]{LevSar:SLDO}
(see also \cite{DymJac:PDECEIP} and \cite{Kre:CATPOUC}).

We note that in \cite[Chap. 12]{LevSar:SLDO} the determination of
a potential $Q$ with prescribed spectral function $\sigma$ 
is incomplete. The self-adjointness of $Q$ is not proved.

\bigskip
The paper is organized as follows. In Section \seceins we present some 
auxiliary
results. In particular  we prove the self-adjointness of the operator $L_H$.

In Section \seczwei we introduce the generalized Fourier transform

$$({\mathcal F}_{H,Q}f)(\lambda ):=\int^{\infty }_0 Y(x,\lambda )^*f(x)dx
   $$
(see \myref{G10-2.19})
for $f\in L^2_{\rm comp}(\R_+,\C^{2n})$ and establish the existence of an 
$n\times n$ matrix
(spectral) measure $\sigma $ such that the Parseval equality
\begin{equation}
(f,g)_{L^2(\R_+,\C^{2n})}=({\mathcal F}_{H,Q}f, \
{\mathcal F}_{H,Q}g)_{L^2_{\sigma }(\R)}, 
  \tag{\ref{G10-0.4}'}
\end{equation}
which is equivalent to \myref{G10-0.4}, holds.
In the proof we follow Krein's method of directing
functionals \cite{Kre:GMDHPNEP}, \cite{Kre:HODF}. 
Moreover we show that ${\mathcal F}_{H,Q}$ is a unitary transformation from
$L^2(\R_+,\C^{2n})$ onto $L^2_{\sigma }(\R)$ which diagonalizes
the operator $L_H$. Namely, 
${\mathcal F}_{H,Q}L_H{\mathcal F}^{-1}_{H,Q}=\Lambda$ where 
$\Lambda:L_\sigma^2(\R)\to L^2_\sigma(\R)$
denotes the multiplication operator by the function $\gl\mapsto \gl$.
Similar results (with similar proofs) hold for Sturm-Liouville operators 
as well as for higher
order differential operators. 

In Section \secdrei we introduce (under the
additional assumptions on $B$)
a triangular transformation operator $I+K$ 
and present a sketch of proof of the representation
$Y(\cdot,\lambda )=((I+K)e_0)(\cdot,\lambda )$ 
where $e_0(x,\lambda )$ is the solution
of \myref{G10-0.3} with $Q=0$. 
Then we derive the linear Gelfand-Levitan equation
\begin{equation}
F(x,t)+K(x,t)+\int^x_0K(x,s)F(s,t)ds=0,\quad  x>t,
\label{G10-0.5}
\end{equation}
with $F(x,t)$ defined by \myref{G2-1.27}.
$F$ is the analog of the so-called transition
function (cf. \cite{Kre:TFODSOBVP}). 
We present two proofs of \myref{G10-0.5}. The proof after Theorem
\plref{S1-1.6} is close to the proofs in \cite{GelLev:DDESF} and
\cite[Chap. 12]{\Lev}. The second one is relatively short.
It is based on simple identities for kernels
of Volterra operators (see \myref{G4-3.3}--\myref{G4-3.9}).
In Proposition \plref{S1-1.5} 
we derive two representations \myref{G4-3.23} and \myref{G2-1.27}
for 
$F(x,t)$ which easily imply \myref{G10-0.5}. 
In other words, this proof derives
the linear equation \myref{G10-0.5} directly from the nonlinear
Gelfand-Levitan equation \myref{G4-3.23}. 
This proof seems to be new and is
essential in the sequel.

Furthermore, in Section \secvier 
we solve the inverse problem (Theorem \plref{S2-4.1}). Namely,
starting with the transition matrix function $F(x,t)$ of the form 
(\ref{G2-2.30}') we prove
the existence of the unique solution $K(x,t)$ of \myref{G10-0.5}. 
Conversely, starting
with $K(x,t)$ we determine the matrix potential $Q(x)=iBK(x,x)-iK(x,x)B$ and
we prove that $Y(\cdot,\lambda ):=((I+K)e_0)(\cdot,\lambda )$ satisfies
the initial value problem \myref{G10-0.3}.


We present several criteria for the prerequisites of Theorem \plref{S2-4.1}
to hold.

Finally, in Section \secfuenf we present some generalizations and improvements
of the main result. The degenerate Gelfand-Levitan equation is also considered
here. We point out that we have obtained a sufficient condition
for an increasing matrix function $\sigma$ to be the spectral function
of the operator $L$. In the special case that $B=(\gl_1 I_n,-\gl_2 I_n)$ 
(or more generally for the class $(T_B)$, cf. Section \secdrei)
our conditions are also necessary. Finally, we prove the existence
of $2n\times 2n$ systems with purely absolute continuous, purely singular
continuous, and purely discrete spectrum of
any given multiplicity $p$, $1\le p\le n$.

In conclusion we mention some recent publications close to our work. 
D. Alpay and I. Gohberg \cite{AlpGoh:ISPDORSMF},\cite{AlpGoh:PARW}
have constructed some explicit formulas
for the  matrix potential of a Dirac system \myref{G1-1.1}
from the rational spectral function.
Their approach is based on the results of minimal factorizations and
realizations of matrix functions \cite{BarGohKaa:MFMOF}. 

A new approach to inverse spectral
problems  for one-dimensional Schr\"odinger operators with partial information
on the potential as well as to different kinds of uniqueness problems on the
half-line has been recently proposed by F. Gesztesy and B. Simon 
(see \cite{GesSim:UTISTODSO}, \cite{GesSim:ISAPIPI} and references therein).
Furthermore, we mention the recent paper F. Gesztesy and H. Holden
\cite{GesHol:TFSTO} on trace formulas for Schr\"odinger-type operators.

The results of this paper have been announced in \cite{LesMal:ISPSHL},
a preliminary version of this paper has been published 
in \cite{LesMal:ISPFOSHL}.

\bigskip
\centerline{\sc Acknowledgements} 
\medskip

The first named author gratefully acknowledges the hospitality
and financial support of the Erwin--Schr\"odinger 
Institute, Vienna, where part of this work was completed.  
Furthermore, the first named author was supported through
the Gerhard Hess Program and the Sonderforschungsbereich 288
of Deutsche Forschungsgemeinschaft.

The second named author gratefully acknowledges the hospitality
   and financial support of the Humboldt University, Berlin,
   where the part of this work was done.


\section{Preliminaries}\label{sec1}

We consider again the operator \myref{G1-1.1} from the introduction.
In the sequel for a vector $v\in\C^{2n}$ the vectors $v_1,v_2\in\C^n$ will
denote the first resp. last $n$ components of $v$.
In this paper scalar products will be antilinear in the first and linear
in the second argument. This is necessary since we will
be dealing with vector measures (see \myref{G2-2.18} below).

$L$ is a formally self--adjoint operator acting on $H^1_{\rm comp}((0,\infty),\C^{2n})
\subset L^2(\R_+,\C^{2n})$. We denote by $L^*$ the adjoint of $L$ in
$L^2(\R_+,\C^{2n})$.
To obtain self--adjoint extensions
we impose boundary conditions of the form
\begin{equation}
   H_2 f_2(0)=H_1 f_1(0).
  \label{G1-1.2a}
\end{equation}
Here, $H_1, H_2\in {\rm M}(n,\C)$ and $f_1(0), f_2(0)\in \C^n$ denote
the first $n$ resp. last $n$ components of $f(0),$
where $f\in H^1_{\rm comp}(\R_+,\C^{2n})$.

\begin{prop}\label{S1-1}
  Let $L_{H_1,H_2}$ be the operator $L^*$ restricted to the domain
$$\cd(L_{H_1,H_2}):=\{ f\in\cd(L^*)\,|\, H_2f_2(0)=H_1f_1(0)\}.$$
Then the operator $L_{H_1,H_2}$ is self--adjoint iff the matrices
$H_1, H_2$ are invertible and $B_1=H^*B_2H$, where $H:=H_2^{-1}H_1$.
\end{prop}

Consequently we have $L_{H_1,H_2}=L_{H_2^{-1}H_1,I}=L_{H,I}$. From now on we
will denote $L_{H,I}$ by $L_H$ and we will write the boundary condition
always in the form
\begin{equation}
    f_2(0)=Hf_1(0).
    \label{G1-1.2}
\end{equation}
\begin{proof} Since $Q$ is continuous we have $\cd(L^*)\subset H_{\rm
loc}^1(\R_+,\C^{2n})$. Now choose a sequence of functions
$\chi_m\in\cinfz{\R}$ with the following properties:
\begin{itemize}
\item[(i)] $\chi_m|(-\infty,m]=1,$
\item[(ii)] $0\le \chi_m\le 1,$
\item[(iii)] $|\chi_m'|\le \frac 1m$.
\end{itemize}
\newcommand{\LtwoCn}{L^2(\R_+,\C^{2n})}
If $f\in\cd(L^*)$ then $\chi_m f\to f$ in $\LtwoCn$ and
\begin{equation}
    L\chi_m f=B \frac 1i \chi_m' f+\chi_m Lf\to Lf
    \label{G1-1.3}
\end{equation}
in $\LtwoCn$. Thus $\chi_m f \to f$ in $\cd(L^*)$.

For $f,g\in\cd(L^*)$ we then find
\begin{equation}
   \begin{array}{cl}
   &(L^*f,g)-(f,L^*g)= \lim_{k\to\infty}\lim_{l\to\infty}
      (L^*\chi_kf,\chi_lg)-(\chi_kf,L^*\chi_lg)\\[0.5em]
     = &-i  \lim_{k\to\infty}\lim_{l\to\infty}
        <B \chi_k f(0),\chi_lg(0)>_{\C^{2n}}
    = -i <Bf(0),g(0)>_{\C^{2n}}.
   \end{array}
    \label{G1-1.4}
\end{equation}

Hence, $g\in \cd(L_{H_1,H_2}^*)$ iff for all $f\in\cd(L_{H_1,H_2})$
\begin{equation}
    0=<Bf(0),g(0)>_{\C^{2n}}.
    \label{G1-1.5}
\end{equation}

This shows that any self--adjoint
extension of $L$ is given by a Lagrangian subspace $V$ of the 
symplectic vector space $\C^{2n}$ with symplectic form
$$\omega(v,w):=<Bv,w>=<B_1v_1,w_1>-<B_2v_2,w_2>.$$
Lagrangian means that $\dim V=n$ and $\omega|V=0$.
The domain of such an extension then is
$$\{ f\in\cd(L^*)\,|\, f(0)\in V\}.$$

Now let $V$ be a Lagrangian subspace of $\C^{2n}=\C^n_+\oplus\C^n_-$. We denote
by $\pi_1, \pi_2$ the orthogonal projections onto the first resp. second
factor. Since the symplectic form $\omega$ is positive resp. negative
definite on $\ker \pi_1$ resp. $\ker \pi_2$ and since $\dim V=n$
the maps $\pi_1, \pi_2$ restricted to $V$ are isomorphisms
\begin{equation}
    \tilde\pi_1:V\to \C^n_+,\quad \tilde\pi_2:V\to \C^n_-.
\end{equation}

Hence $V=\{(x,\tilde\pi_2\circ\tilde\pi_1^{-1}x)\,|\, x\in\C_+^n\}.$
Put $H:=\tilde\pi_2\circ\tilde\pi_1^{-1}$. Then $\go|V=0$ immediately
implies 
$B_1=H^*B_2H.$
This proves the proposition.\end{proof}

\begin{remark}\label{S1-2} 1. The previous proposition shows that 
the deficiency indices $n_\pm(L)$ are equal to $n$, i.e.
$n_\pm(L)=n$. This means that at infinity we do not have
to impose a boundary condition. Thus infinity is always in the
'limit point case', which essentially distinguishes first order
systems from Sturm--Liouville operators and higher order
differential operators (\cite{\Nai,\Cod}).

2. For scalar Dirac systems $(n=B_1=B_2=1)$ another proof of Proposition
\plref{S1-1} has been obtained earlier by B.M. Levitan 
\cite[Theorem 8.6.1]{\Lev}.

The present proof is adapted from the standard proof of the essential
self--adjointness of Dirac operators on complete manifolds
(see e.g. \cite[Theorem II.5.7]{LawMic:SG}).

3. At the same time
as our preprint \cite{LesMal:ISPFOSHL} the paper Sakhnovich
\cite{Sak:DIFOSDE} appeared.
Following Levitan's method he
obtained some sufficient conditions for a canonical system to be
selfadjoint. 
This is a system
\begin{equation}
    J \frac{dy(x,\gl)}{dx}=i H(x) y(x,\gl), \quad
    J=\begin{pmatrix} 0 & I_n\\
                   - I_n& 0 \end{pmatrix},
\label{revision-G3.22}
\end{equation}
where $H(x)$ is a continuous nonnegative $2n\times 2n$ matrix function.
The method of proof of Proposition
\plref{S1-1} can be extended to
arbitrary first order systems, in particular to
generalize the recent result from
\cite{Sak:DIFOSDE} for canonical systems.
Details will be given in a subsequent publication.


4. Another proof of the previous proposition could be given
using the uniqueness of the solution of the Goursat problem
for the hyperbolic system $\frac{du}{dt}=\pm i L_H^*u$ in
$\R_+^2$. This method (see \cite{Ber:EESO})
was also used to prove the essential self--adjointness
of all powers of the Dirac operator on a complete manifold
(cf. \cite{Chernoff}). 
Sakhnovich's result \cite{Sak:DIFOSDE} mentioned before
also follows from the hyperbolic system method.

For the problem considered here we prefered to
present an elementary direct proof.
\end{remark}

From now on we will assume
\begin{equation}
  B_1=H^*B_2H.
  \label{G1-1.8}
\end{equation}
Note that this implies that $H$ is invertible.

We first discuss in some detail the case $Q=0$.
Let $A\in {\rm M}(n,\C)$ be a positive definite matrix. Then we put
for $f\in L^2(\R,\C^n)$
\begin{equation}
   \cf_Af(\gl):= \int_\R e^{-iA^{-1}x\gl}f(x) dx.
   \label{G1-1.9}
\end{equation}
Then we have for $f,g\in L^2(\R,\C^n)$ the Parseval equality
\begin{equation}
   (f,g)=\frac{1}{2\pi}\int_\R (\cf_Af)(\gl)^* A^{-1}(\cf_Ag)(\gl) d\gl.
   \label{G1-1.10}
\end{equation}
To prove \myref{G1-1.10} we may assume $A$ to be diagonal, i.e.
$A={\rm diag}(a_1,\ldots,a_n)$, because if $A=U\tilde A U^*$ with a unitary
matrix $U$ then $(\cf_A f)(\gl)=U(\cf_{\tilde A}U^* f)(\gl)$. 
Now $\cf_Af(\gl)=(\cf
f_j(\gl/a_j))_{j=1,\ldots,n}$ and \myref{G1-1.10} follows easily
from the Parseval equality for the Fourier transform.

Now let 
\begin{equation}
e_0(x,\gl):= e^{i\gl B^{-1} x}{I\choose H}={e^{i\gl B_1^{-1}x}\choose
   e^{-i\gl B_2^{-1} x} H}
   \label{G1-1.11}
\end{equation}
and put
\begin{equation}
   \cf_{H,0}f(\gl):=\int_0^\infty e_0(x,\gl)^* f(x) dx= \cf_{B_1} \tilde f_1(\gl)
+H^* \cf_{B_2}\tilde f_2(-\gl),
   \label{G1-1.12}
\end{equation}
where $\tilde f_j$ denotes the extension by $0$ of $f_j$ to $\R$.

If $f,g\in L^2(\R_+,\C^{2n})$ then the integrals
\begin{equation}\begin{array}{l}
   \DST 
\int_\R \cf_{B_1}\tilde f_1(\gl)^* B_1^{-1} H^*\cf_{B_2} \tilde g_2(-\gl)
   d\gl,\\[1em]
\DST \int_\R \cf_{B_2}\tilde f_2(-\gl)^* H B_1^{-1} \cf_{B_1} \tilde g_1(\gl)
   d\gl
   \end{array}\label{G1-1.13}\end{equation}
are sums of scalar products of the form
\begin{equation}
\int_\R \overline{\cf\tilde\varphi(-\gl)} \cf \tilde \psi(\gl)d\gl,
  \label{G1-1.14}
\end{equation}
where $\varphi,\psi\in L^2(\R_+)$. These scalar products vanish
and hence
we end up with the Parseval equality in the case of $Q=0$
\begin{align}
   \frac{1}{2\pi}\int_\R &\cf_{H,0} f(\gl)^* B_1^{-1} \cf_{H,0} g(\gl) d\gl
    \nonumber\\
   &=\frac{1}{2\pi}\int_\R \cf_{B_1}\tilde f_1(\gl)^* B_1^{-1} 
     \cf_{B_1}\tilde  g_1(\gl) d\gl+\nonumber\\
   &\qquad+\frac{1}{2\pi}\int_\R \cf_{B_2}\tilde f_2(-\gl)^* HB_1^{-1}H^* 
     \cf_{B_2}\tilde  g_2(-\gl) d\gl\nonumber\\
   &= (f,g),\label{G1-1.15}
\end{align}
in view of \myref{G1-1.8} and \myref{G1-1.10}.

\section{The spectral measure}\label{sec2}

In this section we prove the existence of a spectral measure function
for the self--adjoint operator $L_H$ based on Krein's method of directing
functionals \cite{Kre:GMDHPNEP}, \cite{Kre:HODF}. 
For the convenience of the reader we recall Krein's result.

\begin{dfn}[{\rm \cite{Kre:GMDHPNEP}, \cite{Kre:HODF}}]
\label{revision-S3.1} 
Let $A$ be a symmetric operator in a separable Hilbert
space $H$ and let $E$ be a dense linear subspace of
$H$ containing $\cd (A).$

The system $\{\Phi _j\}^p_1$ of linear functionals defined on $E$ and
depending on $\lambda \in \R$ is called a directing system of functionals 
for $A$ in $E$ if the following three conditions are fulfilled:
\begin{enumerate}
\renewcommand{\labelenumi}{{\rm \arabic{enumi}.}}
\item $\Phi_j(f;\lambda ), j=1,...,p$, is an analytic function of
$\lambda \in \R$, for each $f\in E;$
\item
the functionals $\Phi_j(\cdot;\lambda_0)$ are linearly independent for some
$\lambda _0\in \R;$
\item for each $f_0\in E$ and $\lambda_0\in \R$ the equation
$Ag-\lambda _0g=f_0$ has a solution in $E$ if and only if
$$
\Phi_j(f_0;\lambda _0)=0\quad\mbox{\rm for all } 1\le j\le p.
$$
\end{enumerate}
\end{dfn}

\begin{theorem}[{\rm \cite{Kre:GMDHPNEP}, \cite{Kre:HODF}}]
\label{revision-S3.2}
Let $A$ be a symmetric operator in $H$ with
$\cd(A)\subset E \subset H$ which has a directing system of
functionals $\{\Phi_j(\cdot;\lambda )\}^p_1$ in $E$. 
Then

{\rm 1.} there exists an increasing $p\times p$ matrix function
$\sigma(\lambda )=(\sigma_{jk}(\lambda))^p_{j,k=1}$ 
such that the equality
$$
(g,f)=\sum^p_{j,k=1}\int _{\R}\overline{\Phi _j(g;\lambda )}
\Phi_k(f;\lambda )d\sigma_{jk}(\lambda)
$$
holds for each $f,g\in E$.

{\rm 2.} If $\sigma$ is normalized by requiring it to be right--continuous
with $\sigma(0)=0$ then it is unique if and only if 
$n_+(A)=n_-(A)$, where 
$n_\pm(A):=\dim\ker(A^*\mp i)$ denote the deficiency indices of $A$.
\end{theorem}

\begin{dfn}\label{revision-S3.3}
Let $\sigma(\gl)=(\sigma_{ij}(\gl))_{i,j=1}^n$ be an increasing
$n\times n$ matrix function. On the space $C_0(\R,\C^n)$ of continuous
$\C^n$--valued functions with compact support we introduce the scalar product
\begin{equation}
  (f,g)_{L^2_\sigma}:=\int_\R f(\gl)^* d\sigma(\gl) g(\gl)
    :=\sum_{i,j=1}^n \int_\R
     \overline{f_i(\gl)}g_j(\gl)d\sigma_{ij}(\gl).\label{G2-2.18}
 \end{equation}
We denote by $L^2_\sigma(\R)$ (cf. \cite{\Nai}) the
Hilbert space completion of this space.
\end{dfn}

\begin{remark}
From now on we will consider -- without saying this explicitly --
only right--continuous
$n\times n$ matrix functions which map $0$ to the $0$--matrix.
Such a function $\sigma$ is determined
by its corresponding matrix measure $d\sigma$.
\end{remark}

We turn to general $Q$. For future reference we state the boundary
value problem for $L$:
\begin{equation}
    Lf=\gl f,\quad f_2(0)=H f_1(0),\quad\mbox{\rm where}\quad
    B_1=H^*B_2H.
   \label{G3-2.1}
\end{equation}

\newcommand{\fouHQ}{\cf_{H,Q}}
\begin{prop}\label{S1-1.2}
Let $Y:\R_+\times \C\to {\rm M}(2n\times n,\C)$ be
the unique solution of the initial value problem
\begin{equation}
   LY(x,\gl)=\gl Y(x,\gl),\quad Y(0,\gl)={I\choose H}.
   \label{G1-1.16}
\end{equation}
Then:

{\rm 1.} There exists an increasing
$n\times n$ matrix function
$\sigma(\gl), \gl\in\R,$ (spectral function) such that the map
\begin{equation}
\fouHQ:L^2_{\rm comp}(\R_+,\C^{2n})\ni
f\mapsto (\fouHQ f)(\gl):=F(\gl):=\int_0^\infty Y(x,\gl)^* f(x)
dx
    \label{G10-2.19}
\end{equation}
extends by continuity to an isometric transformation from $L^2(\R_+,\C^{2n})$
into the space $L^2_\sigma(\R)$, i.e. for $f,g\in L^2(\R_+,\C^{2n})$ we have
the Parseval equation
\begin{equation}
     \int_0^\infty f^*(t)g(t) dt = \int_\R F^*(\gl)d\sigma(\gl)G(\gl)
     \label{G2-2.19}
\end{equation}
with $F,G$ being the $\fouHQ$--transforms of $f,g$.

{\rm 2.} If $\sigma$ is normalized by requiring it to be right--continuous
with $\sigma(0)=0$ then it is unique.
\end{prop}

\begin{proof} 1. Let $b\in \R_+$ be a fixed point and 
let $L_b$ be the operator $L^*$
restricted to the domain
\begin{equation}
{\mathcal D}(L_b)=\{f\in H^1([0,b], \C^{2n})\,|\, f_2(0)=Hf_1(0), 
f_1(b)=f_2(b)=0\}.
\end{equation}
It is clear that $L_b$ is a symmetric operator and
$L^*_b=(L_b)^*$ is a restriction of $L^*$ to the domain
\begin{equation}
{\mathcal D}(L_b^*)=\{f\in H^1([0,b], \C^{2n})\,|\, f_2(0)=Hf_1(0) \}.
\end{equation}

We consider ${\mathcal D}(L_b)$ as a subset of $H^1({\R}_+, {\C}^{2n})$
identifying each function $f\in {\mathcal D}(L_b)$ with its continuation by
zero to ${\R}_+$.

Since $L_b$ is a regular differential operator on a finite interval, each
$\lambda \in {\C}$ is a regular type point for $L_b$, i.e.
$\|(L_b-\lambda )f\|\ge\epsilon \|f\|$ for all $f\in {\mathcal D}(L_b)$ with some
$\epsilon>0$. In particular, $L_b-\lambda$ has closed range.
Hence, for a fixed $\lambda \in {\R}$ and
$f\in L^2([0,b];{\C}^{2n})$ the equation
\begin{equation}
Lg-\lambda g=f,\quad \lambda \in {\R},
\end{equation}
has a solution $g\in L^2([0,b],{\C}^{2n})$ if and only if $f$ is
orthogonal to the kernel $\ker (L_b^*-\lambda )$, that is if
$\int ^b_0Y^*(x,\lambda )f(x)dx=0.$

Denoting by $Y_i$ the $i$-th column of $Y$, on rewrites the
last equation as
\begin{equation}\begin{array}{rcl}
\Phi_i(f;\lambda )&:=&\DST\int ^b_0<Y_i(x,\lambda ),f(x)>dx\\[1em]
      &=&\DST\sum ^{2n}_{j=1}\int ^b_0
\overline {Y_{ji}(x,\lambda )}f_j(x)dx=0,\quad 1\le i\le n.
\label{revision-G3.9}
		\end{array}
\end{equation}
It is clear that the functionals $\Phi_i$
on $L^2_{\rm comp}({\R}_+,{\C}^{2n})$,
defined by the left-hand
side of \myref{revision-G3.9},
are linearly independent and holomorphic in
$\lambda \in {\R}$. Thus the conditions 1. and 2.
of Definition \plref{revision-S3.1} are
satisfied. Since $E :=L^2_{\rm comp}({\R}_+,{\C}^{2n})$ is dense in
$L^2({\R}_+,{\C}^{2n})$ the functionals
$\Phi _i(f,\lambda )$ thus form a directing
system of functionals for the operator
$A:=L_H\restriction\cd(L_H)\cap L^2_{\rm comp}({\R}_+,{\C}^{2n})$. 
By Krein's Theorem \plref{revision-S3.2} there
exists $\sigma (\lambda )$ such that  \myref{G2-2.19} holds for arbitrary
$f,g\in L^2_{\rm comp}({\R}_+,{\C}^{2n})$.

2. In view of Proposition \plref{S1-1} the operator
$A=L_H\restriction\cd(L_H)\cap L^2_{\rm comp}({\R}_+,{\C}^{2n})$
is essentially selfadjoint and
consequently $n_+(A)=n_-(A)=0$. 
Thus the uniqueness of $\sigma (\lambda )$
follows from the assertion 2. of Krein's theorem \plref{revision-S3.2}.
\end{proof}

\begin{remark}
1. The Parseval identity may be symbolically rewritten as
\begin{equation}
    \int_{\R} Y(x,\gl)d\sigma(\gl) Y(t,\gl)^*=\delta(x-t) I_{2n}.
    \tag{\ref{G2-2.19}'}
\end{equation}
To obtain (\ref{G2-2.19}') from \myref{G2-2.19} it suffices to
set in \myref{G2-2.19}
  $f(\xi )=\delta _x(\xi )\otimes e_i, \ g(\xi )=\delta _t(\xi )\otimes e_j,
1\le i,j\le 2n$
and to note that
 $(\cf_{H,Q}f)(\gl )=Y(x,\gl )^*e_i, \ (\cf_{H,Q}g)(\gl )=Y(t,\gl )^*e_j$.

2. Another proof of Proposition \plref{S1-1.2} based on the
approximation method proposed independently by B.M. Levitan
\cite[Chap. 8]{LevSar:SLDO} and N. Levinson \cite[Chap. 9]{CodLev:TODE}
was given in the preliminary version of this paper \cite{LesMal:ISPFOSHL}.
\end{remark}

For convenience we denote the extension of $\fouHQ$ to
$L^2({\R}_+,{\C}^{2n})$ by the same letter.
Next we prove the surjectivity of $\fouHQ$.

\begin{theorem}\label{revision-S3.7}
Under the assumptions of Proposition \plref{S1-1.2} the mapping
$\fouHQ$ is surjective, that is ${\mathcal F}_{H,Q}$ maps
$L^2({\R}_+,{\C}^{2n})$ onto $L^2_{\sigma }({\R})$.
\end{theorem}

\begin{proof}
So far we have proved that $\fouHQ:L^2(\R_+,\C^{2n})\longrightarrow
L^2_\sigma(\R)$ is an isometry. To prove surjectivity we mimick
the proof of \cite[Sec. 9.3]{\Cod} for second order operators.

Note first that for $f\in \cd(L_H)$
we have
\begin{equation}
    (\fouHQ L_H f)(\gl)=-Y(0,\gl)^* \frac 1i Bf(0)+\gl (\fouHQ f)(\gl)
       = \gl (\fouHQ f)(\gl)\label{G4-2.22}
\end{equation}
since in view of \myref{G3-2.1} and \myref{G1-1.16}
$Y(0,\gl)^*B f(0)=0.$

For $f\in L^2_{\rm comp}(\R_+,\C^{2n})\cap \cd(L_H)$
formula \myref{G4-2.22} follows from integration by parts. For
arbitrary $f\in\cd(L_H)$ it follows from the fact
that $L^2_{\rm comp}(\R_+,\C^{2n})\cap \cd(L_H)$
is a core for $L_H$. The latter follows from the proof of Proposition
\plref{S1-1}.

Next we construct the adjoint of $\fouHQ$:
we put for $g\in L^2_{\sigma,{\rm comp}}(\R)$
\begin{equation}
   (\cg_H g)(x):=\int_\R Y(x,\gl) d\sigma(\gl) g(\gl).
   \label{G4-2.23}
\end{equation}
Then for $f\in L^2_{\rm comp}(\R_+,\C^{2n})$
\begin{eqnarray}
   (\cg_H g,f)_{L^2(\R_+,\C^{2n})} &=&  \int_0^\infty (\cg_H g)(x)^* f(x) dx
      \nonumber\\
   &=& \int_0^\infty \int_\R g(\gl)^* d\sigma(\gl) Y(x,\gl)^* f(x) dx
       \label{G4-2.24}\\
    &=& (g,\fouHQ f)_{L^2_\sigma(\R)}.\nonumber
\end{eqnarray}
From the estimate
\begin{equation}
      \big|  (\cg_H g,f)_{L^2(\R_+,\C^{2n})}\big|=
    \big|(g,\fouHQ f)_{L^2_\sigma(\R)}\big|\le
    \|g\|_{L^2_\sigma(\R)}\|f\|_{L^2(\R_+,\C^{2n})}
    \label{G4-2.25}
\end{equation}
we infer that $\cg_H$ extends by continuity for $L^2_\sigma(\R)$.
Moreover, it equals the adjoint of $\fouHQ$, i.e.
\begin{equation}
    \cg_H=\fouHQ^*.
   \label{G4-2.26}
\end{equation}
Since $\fouHQ$ is an isometry it remains to prove injectivity of
$\fouHQ^*$.

It follows from \myref{G4-2.22} and Proposition \plref{S1-1}
\begin{equation}
\fouHQ(L_H-\zeta)^{-1}=(\Lambda-\zeta)^{-1}\fouHQ
  \label{G6-2.27}
\end{equation}
for $\zeta=\nu+i\varepsilon \in \C\setminus \R$.
Here $\Lambda:L^2_\sigma(\R)\to L^2_\sigma(\R), (\Lambda g)(\gl):=\gl g(\gl)$
denotes the operator of multiplication by $\lambda$.
Therefore $(L_H-\bar\zeta)^{-1}\fouHQ^*=
\fouHQ^*(\Lambda-\bar\zeta)^{-1}$ and hence we have the implication:
\begin{equation}
   \Phi\in \ker \fouHQ^*\Longrightarrow
   (\Lambda -\zeta)^{-1}\Phi \in \ker \fouHQ^*
   \quad\mbox{\rm for all}\; \zeta\in\C\setminus\R.
   \label{G4-2.27}
\end{equation}

We put
\begin{equation}
       \widetilde Y(x,\gl):=\int_0^x Y(t,\gl)dt.
   \label{G4-2.28}
\end{equation}
Note that the $i$--th row $\widetilde Y_i(x,\gl)$ of $\widetilde Y$ satisfies
\begin{equation}
       \widetilde Y_i(x,\gl)^*=\int_0^x Y(t,\gl)^* e_i dt
       = (\fouHQ (1_{[0,x]}\otimes e_i))(\gl),
   \label{G4-2.29}
\end{equation}
where $e_i$ denotes the $i$--th unit vector in $\C^{2n}$.

Now let $\Phi\in\ker \fouHQ^*$.

In particular $\widetilde Y_i(x,\cdot)^*\in L^2_\sigma(\R)$ and thus
in view of \myref{G4-2.24} and \myref{G4-2.27} we have for
$x\ge 0, \eps>0, \nu\in\R$
\begin{equation}
   0=(1_{[0,x]}\otimes e_i,\fouHQ^* \frac{\eps}{(\Lambda-\nu)^2+\eps^2}\Phi)
    =(\widetilde Y_i^*,  \frac{\eps}{(\Lambda-\nu)^2+\eps^2}\Phi),
\end{equation}
thus
\begin{equation}
     0=\int_\R \widetilde Y(x,\gl) \frac{\eps}{(\gl-\nu)^2+\eps^2} d\sigma(\gl) \Phi(\gl).
\end{equation}
Since $\widetilde
Y(x,\cdot)d\sigma(\cdot) \Phi(\cdot)$ is $L^1$ the dominated
convergence theorem implies for $\ga,\gb\in\R$
\begin{eqnarray}
    0&=&   \lim_{\eps\to 0} \int_\ga^\gb \int_\R \widetilde Y(x,\gl)
         \frac{\eps}{(\gl-\nu)^2+\eps^2}d\sigma(\gl) \Phi(\gl)d\nu\nonumber\\
   &=& \int_\R \widetilde Y(x,\gl) \lim_{\eps\to 0}\int_\ga^\gb
       \frac{\eps}{(\gl-\nu)^2+\eps^2}d\nu d\sigma(\gl) \Phi(\gl)\\
   &=& \pi \int_\ga^\gb \widetilde Y(x,\gl)d\sigma(\gl) \Phi(\gl).\nonumber
\end{eqnarray}
Differentiating by $x$ and putting $x=0$ yields for $\ga,\gb\in\R$
 $$   0= \int_\ga^\gb Y(0,\gl) d\sigma(\gl) \Phi(\gl).$$
Since $Y(0,\gl)={I\choose H}$ and $H$ is invertible we have
 $$    \int_\ga^\gb d\sigma(\gl) \Phi(\gl)=0$$
for all $\ga,\gb\in \R$. This implies $\Phi=0$ in $L^2_\sigma(\R)$. 

\end{proof}

\begin{remark}\label{revision-S2.2}
1.
We note that another proof of the uniqueness of the spectral function
in Proposition \plref{S1-1.2} can be given using Theorem \plref{revision-S3.7}.
To prove the uniqueness statement we assume we had another increasing right
continuous $n\times n$ matrix function $\varrho, \varrho(0)=0$, such that
$\fouHQ$ is a unitary transformation from $L^2(\R_+,\C^{2n})$ onto
$L^2_\varrho(\R)$. Then in view of \myref{G2-2.19} we have
for all $F,G\in L^2_{\rm comp}(\R,\C^n)$
$$\int_\R F(\gl)^* d\sigma(\gl) G(\gl)=
  \int_\R F(\gl)^* d\varrho(\gl) G(\gl),$$
and hence the two Radon vector measures $d\sigma$ and $d\varrho$
coincide. By the right--continuity and the normalization
$\varrho(0)=\sigma(0)=0$ this implies $\sigma=\varrho$.

2. For $n=1$ Proposition \plref{S1-1.2} follows from
\cite[Theorem 4]{Kre:HODF}. We also note that a generalization of 
Krein's theorem to the case $n>1$ may be obtained by a slight
modification of the proof of Proposition \plref{S1-1.2}.

3. 
In \cite[Chap. 3]{Sak:FPOI}
the existence of the spectral function for a canonical system
\myref{revision-G3.22}
is stated. For nonsingular Hamiltonians
this fact follows from Krein's Theorem 
\plref{revision-S3.2}
in just the same way as Proposition \plref{S1-1.2}. 

We note also that for a singular Hamiltonian similar results may be obtained
by the corresponding generalization of Krein's Theorem 
\plref{revision-S3.2} for linear relations.
\end{remark}

\begin{example} \myref{G1-1.15} shows that in the case $Q=0$ we can choose
for $\sigma$ the function $\sigma_0(\gl):=\frac{1}{2\pi} B_1^{-1} \gl$.

\end{example}

\section{Transformation operator and Gelfand--Levitan equation}
\label{sec3}

{\bf 1.}\quad
We present a special case of \cite[Theorem 7.1]{\Mal},
(see also \cite[Theorem 1.2]{Mal:UQIPFOSFI}).
In the sequel we assume $B$ to be a diagonal matrix, which can
be achieved by conjugating $L$ with an appropriate unitary matrix.

Let
\begin{equation}\begin{split}
       &B=\diag(\gl_1 I_{n_1}, \ldots,\gl_r I_{n_r}),\\
       &n_1=\min\{n_i\,|\,1\le i\le r\},\quad
       n_1+n_2+\ldots+ n_r=2n.\label{mal3.1}
    		\end{split}
\end{equation}
Furthermore, we put
\begin{equation}
     \gO:=\{(x,t)\in \R^2\,|\, 0\le t\le x\}.\label{G4-3.1}
\end{equation}

\begin{theorem}\label{revision-S4.1}
Let $B$ be as in \myref{mal3.1} and
let $Q=(Q_{ij})_{i,j=1}^r:\R_+\longrightarrow
{\rm M}(2n,\C)$ continuous, where $Q_{ij}$ denotes the block--matrix
decomposition with respect to the orthogonal decomposition
$\C^{2n}=\moplus\limits_{i=1}^r \C^{n_i}$.
Moreover, we assume that $Q$ is off--diagonal, i.e.
\begin{equation}
      Q_{ii}=0,\quad i=1,\ldots,r.
   \label{mal3.2}
\end{equation}
Let $Y$ be the solution of the equation \myref{G1-1.16} satisfying
\begin{equation}
  Y(0,\gl)=A:=\operatorname{col}(A_1, \ldots, A_r), 
    \quad A_j \in {\rm M}(n_j\times n_1, \C),
\quad \rank A_j=n_1.
  \label{revision-G4.4}
\end{equation}
Then there exists a continuous function $K:\gO\longrightarrow
{\rm M}(2n,\C)$ such that we have
\begin{equation}
  Y(x,\gl)=Y_0(x,\gl)+\int_0^x K(x,t) Y_0(t,\gl) dt,
       \label{mal3.3}
\end{equation}
where 
\[Y_0(x,\gl)=e^{i\gl B^{-1}x} A\]
is the solution of the equation 
\myref{G3-2.1} with $Q=0$ and
satisfying the same initial conditions \myref{revision-G4.4}.

If $Q\in C^1(\R_+,{\rm M}(2n,\C))$ then $K\in C^1(\gO,{\rm M}(2n,\C))$
and it satisfies
\alpheqn[G4-3.15]
\begin{eqnarray}
&&   B\pl_xK(x,t)+\pl_tK(x,t)B+iQ(x) K(x,t)=0,\label{G4-3.15a}\label{mal3.4a}\\
&& BK(x,x)-K(x,x) B+iQ(x)=0,\label{G4-3.15b}\label{mal3.4b}\\
&&K(x,0)BA=0.\label{G4-3.15c}\label{mal3.4c}
\end{eqnarray}
\reseteqn
If $Q\in C(\R_+,{\rm M}(2n,\C))$ then $K$ is the generalized continuous
solution of \myref{G4-3.15}.

Conversely, if $K$ is a (generalized) solution of \myref{G4-3.15} then
$Y(x,\gl)$ defined by \myref{mal3.3} is the (generalized) solution of the
initial value problem \myref{G1-1.16}. 

\end{theorem}
\begin{proof}[Sketch of proof]
i) Suppose that $K\in  C^1(\R_+,{\rm M}(2n,\C))$ and that formula
\myref{mal3.3} holds. Substituting \myref{mal3.3} into \myref{G1-1.16}
and integrating by parts one obtains
\begin{equation}\begin{split}
 &\big[BK(x,x)-K(x,x) B+iQ(x)\big] Y_0(x,\gl) +
     K(x,0)BY_0(0,\gl)\\
 & +\int_0^x\big[B\pl_xK(x,t)+\pl_tK(x,t)B+i Q(x) K(x,t)\big]Y_0(t,\gl)dt=0.
		\end{split}
    \label{G4-3.16}
\end{equation}
Since $Y_0(0,\gl)=A$ does not depend on $\gl$ one concludes
from \myref{G4-3.16} and the Riemann--Lebesgue Lemma that \myref{G4-3.16}
is equivalent to \myref{G4-3.15}. Thus in this
case the representation \myref{mal3.3} is equivalent to the solvability
of the problem \myref{G4-3.15}. 

\newcommand{\lam}{\lambda}

ii) Next we prove the existence of a (not unique) solution
of the problem \myref{mal3.4a}--\myref{mal3.4b}.
Let $R(x,t)$ be one of them. Using the block--matrix
representation $R(x,t)=(R_{ij}(x,t))^r_{i,j=1}$ we rewrite the problem
\myref{mal3.4a}--\myref{mal3.4b} as
\begin{align}
 \lam_i \pl_xR_{ij}(x,t)+\lam _j\pl_tR_{ij}(x,t)
   &=-\sqrt{-1}\sum ^r_{p=1}Q_{ip}(x)R_{pj}(x,t),
  \quad  1\le i,j\le r,
    \label{mal3.6}\\
R_{ij}(x,x)&=
  -\sqrt{-1}(\lam _i-\lam _j)^{-1} \  Q_{ij}(x), \quad  1\le i\not= j\le r.
   \label{mal3.7}
\end{align}
It is clear that the system \myref{mal3.6}
is hyperbolic with real characteristics
$l_{ij}:x=k_{ij}t+c (k_{ij}=\lam_j \lam ^{-1}_i)$. Thus, in
$\Omega=\{0\le t\le x<\infty\}$ we
have the incomplete characteristic Cauchy problem
\myref{mal3.6}, \myref{mal3.7} with
$(2n)^2-n_1^2-\ldots-n_r^2$ scalar conditions \myref{mal3.7}.
Fixing $x_0\in \R_+$ and setting
$$
k_{\min}=
  \max\{k_{ij}\,|\,\ k_{ij}\in (0,1),\ 1\le i,j\le r\},
 \quad k_{\max}=k_{\min}^{-1},
$$
we consider the triangle $\triangle_{ABC}$ confined by the lines
$AB:\ x=t,\ AC:\ x-x_0=k_{\min}t,\ BC:\ x-x_0=k_{\max}t.$ We preserve
the notation $Q(x)$ for a continuous extension to
$\R$ of the function $Q(x)$ with the same norm.
Furthermore, we denote by $a$ and $b$ the abscissas
of the points $A$ and $B$ respectively.
Now we impose the following $n_1^2+\ldots+n_r^2$ conditions on the
characteristic line $AC:$
\begin{equation}
R_{jj}(x,(x-1)k_{\min})=0,
\quad\mbox{for}\quad x\in [a,x_0],  \quad j \in {1,...,r}.
\label{mal3.8}
\end{equation}
Thus, we arrive at the Goursat problem \myref{mal3.6}--\myref{mal3.8}
for the hyperbolic
system \myref{mal3.6} in the triangle $\triangle_{ABC}.$
Integrating the system \myref{mal3.6}
along the characteristics and using \myref{mal3.7},
\myref{mal3.8} one deduces the system of integral equations
\begin{equation}\begin{split}
\lam_i R_{ij}&(x,t)\\
 &=\frac{\lam_i Q_{ij}(\xi_{ij}(x,t))}{\lam_i-\lam_j}-
\sqrt{-1}
\int_{\xi_{ij}(x,t)}^x\sum_{p=1}^rQ_{ip}(\xi)R_{pj}(\xi,(\xi-x)k_{ij}+t)d\xi,
\label{mal3.9}
		\end{split}
\end{equation}
where for brevity it is set 
$\lam_i(\lam_i-\lam_j)^{-1}Q_{ij}(\xi_{ij}(x,t))=0$ for
$i=j$ and
\[
\xi_{ij}(x,t)=
\begin{cases}(\lam_jx-\lam_it)(\lam_j-\lam_i)^{-1},&i\not= j,\\
 a+(1-a)(x-t),&i=j.
\end{cases}
\]
For $Q\in C^1(\R,{\rm M}(n,\C)$ the system \myref{mal3.9}
is equivalent to the Goursat problem
\myref{mal3.6}--\myref{mal3.8}. The solvability (and uniqueness)
of the solution of \myref{mal3.9}
is proved by the method of successive approximations.

For $Q\in C(\R,{\rm M}(n,\C))\setminus C^1(\R,{\rm M}(n,\C)$ we understand
the solution of \myref{mal3.6}--\myref{mal3.8} as a solution of \myref{mal3.9}.

iii) To finish the proof, starting with the solution $R(x,t)$ of the Goursat
problem \myref{mal3.6}--\myref{mal3.8} we introduce a convolution operator
$$
\Phi: f\to \int_0^x\Phi(x-t)f(t)dt
$$
with $\Phi(x)=\diag(\Phi_1(x),\ldots,\Phi_r(x))$ being a block--diagonal
$2n\times 2n$ matrix function, consisting of $n_j\times n_j$ blocks $\Phi_j$
and define the operator $K$ by the equality $I+K=(I+R)(I+\Phi).$
It is clear that $K$ is a Volterra operator with the kernel
\begin{equation}
   K(x,t)=R(x,t)+\Phi(x-t)+\int_0^xR(x,s)\Phi(s-t)ds.
   \label{mal3.10}
\end{equation}
Since the operator $I+R$ intertwines the restrictions $L_0$ and
$-iB\otimes D_0$ of the operators $L$ and $-iB\otimes D$ onto
$\{f\in H^1([0,1],\C^{2n})\,|\, f(0)=0\},$
that is $L_0(I+R)=(I+R)(-iB\otimes D_0),$ so is $I+K.$
This fact amounts to saying that $K(x,t)$ satisfies the problem
\myref{mal3.4a}--\myref{mal3.4b}. To satisfy
the condition \myref{mal3.4c} it suffices
(in view of \myref{mal3.10}) to choose $\Phi(x)$
as the solution of the equation
\begin{equation}
\Phi(x)BA +\int_0^x R(x,s)\Phi(s)BA ds=-R(x,0)BA.
\label{mal3.11}
\end{equation}
Since $\rank A_j=n_1, 2\le j\le r$,  the Volterra
equation \myref{mal3.11} is of the second kind and therefore
has the unique solution
$\Phi\in C([0,\infty),{\rm M}(2n,\C)).$ Thus
$K(x,t)$ is the required solution of \myref{mal3.4a}--\myref{mal3.4c}.
\end{proof}

\begin{cor}\label{revision-S4.2}
Under the assumptions of Theorem \plref{revision-S4.1} let 
$B=(\gl_1 I_n, \gl_2 I_n)$ (that is $r=2.$)

Then there exists a continuous function $K:\gO\longrightarrow
{\rm M}(2n,\C)$ such that we have
\begin{equation}
  Y(x,\gl)=e_0(x,\gl)+\int_0^x K(x,t) e_0(t,\gl) dt,
   \label{revision-G4.13}
\end{equation}
where $e_0(x,\gl)$ was defined in \myref{G1-1.11}.

If $Q\in C^1(\R_+,{\rm M}(2n,\C))$ then $K\in C^1(\gO,{\rm M}(2n,\C))$
and it satisfies
\alpheqn[revision-G4.14]
\begin{eqnarray}
&&   B\pl_xK(x,t)+\pl_tK(x,t)B+iQ(x) K(x,t)=0,\label{revision-G4.14a}\\
&& BK(x,x)-K(x,x) B+iQ(x)=0,\label{revision-G4.14b}\\
&&K(x,0)B{I\choose H}
=0.\label{revision-G4.14c}
\end{eqnarray}
\reseteqn
If $Q\in C(\R_+,{\rm M}(2n,\C))$ then $K$ is the generalized continuous
solution of \myref{G4-3.15}.

Conversely, if $K$ is a (generalized) solution of \myref{revision-G4.14} then
$Y(x,\gl)$ defined by \myref{revision-G4.13}
is the (generalized) solution of the
initial value problem \myref{G1-1.16}. This last statement holds even for
general $B$ of the form \myref{G1-1.1}.
\end{cor}


\medskip\noindent
{\bf 2.}\quad
We continue with some general remarks about Volterra operators:

For any continuous matrix function $K:\gO\longrightarrow {\rm M}(2n,\C)$
we obtain a Volterra operator
\begin{equation}
     Kf(x):=\int_0^x K(x,t) f(t) dt
     \label{G4-3.2}
\end{equation}
acting on $C(\R_+,\C^{2n})$ or $L^2([0,a],\C^{2n})$ for any $a>0$.
By slight abuse of notation we will use the same symbol for the
operator and its kernel.
The set of operators $I+K$ with $K$ being a Volterra operator forms a group.
The operator
\begin{equation}
    R:= (I+K)^{-1}-I
    \label{G4-3.3}
\end{equation}
is again a Volterra operator with continuous kernel $R(x,t), t\le x$.
From the equation
\begin{equation}
     I=(I+R)(I+K)=(I+K)(I+R)
      \label{G4-3.4}
\end{equation}
we deduce
\begin{equation}
          RK=KR=-R-K.\label{G4-3.5}
\end{equation}
Put 
\begin{equation}
      F:=R+R^*+RR^*.\label{G4-3.6}
\end{equation}
The kernel of $F$ obviously is
\begin{equation}
       F(x,t)=\left\{\begin{array}{ll}
        \DST  R(x,t)+\int_0^t R(x,s)R(t,s)^*ds,& x>t,\\[1em]
        \DST  R(t,x)^* +\int_0^x R(x,s) R(t,s)^* ds,& x<t.
        	      \end{array}\right.
        \label{G4-3.7}    
\end{equation}
Furthermore, using \myref{G4-3.5} we conclude
\begin{equation}
       F+K+KF=R+K+R^*+RR^*+KR+KR^*+KRR^*=R^*
       \label{G4-3.8}
\end{equation}
thus we have the "Gelfand--Levitan equation"
\begin{equation}
         F+K-R^*+KF=0.\label{G4-3.9}
\end{equation}

\begin{prop}\label{S4-3.1} Let $K:\gO\longrightarrow {\rm M}(2n,\C)$ be continuous
and let $R:\gO\longrightarrow {\rm M}(2n,\C)$ be the continuous kernel of
the Volterra operator $(I+K)^{-1}-I$. Then the function 
$F:\R_+^2\longrightarrow {\rm M}(2n,\C)$ defined by \myref{G4-3.7}
satisfies the ``Gelfand--Levitan equation''
\begin{eqnarray}
  &&F(x,t)+K(x,t)+\int_0^x K(x,s) F(s,t) ds =0,\quad x>t,\label{G4-3.10}
        \\[0.5em]
  &&F(x,t)-R(t,x)^*+\int_0^x K(x,s) F(s,t) ds =0,\quad x<t.\label{G4-3.11}
\end{eqnarray}
Conversely, if $F_1:\gO\longrightarrow {\rm M}(2n,\C)$ is continuous
and satisfies \myref{G4-3.10} then $F_1=F|\gO$.
\end{prop}
\begin{proof} It only remains to prove the assertion about $F_1$. The difference
$F(x,t)-F_1(x,t)$ satisfies the equation
$$ F(x,t)-F_1(x,t)+\int_0^x K(x,s)[F(s,t)-F_1(s,t)]ds=0,\quad 0\le t\le x.$$
For each fixed $t\in [0,x]$ this is a homogeneous Volterra equation of
the second kind and consequently has only the trivial solution
$F(x,t)-F_1(x,t)=0$.
\end{proof}

We turn back to the system \myref{G3-2.1}. 

\begin{dfn} We say that the system
\myref{G3-2.1} (resp. the operator $L$)
belongs to the class $(T_B)$ if for this system
there exists a transformation operator. 
\end{dfn}

This means that the solution $Y(x,\lambda )$ of the initial value
problem \myref{G1-1.16}
admits a representation \myref{revision-G4.13} with a continuous function
$K:\Omega \to M(2n;\C)$. Corollary \plref{revision-S4.2} says that
the system is of class $(T_B)$ if $B=(\gl_1 I_n, \gl_2 I_n)$. 

It follows easily from Proposition \plref{S1-1.5} below
that for an operator $L_H$
of class $(T_B)$ the transformation operator $I+K$ is unique, i.e. the
representation \myref{G1-1.16} for $Y(x,\lambda )$ is unique.

If the system is of class $(T_B)$ then we 
denote by $K$ the unique Volterra operator with continuous kernel
satisfying \myref{revision-G4.14}. As before $R$ denotes the Volterra operator
defined by $R:=(I+K)^{-1}-I$. 

In particular we have in view of \myref{revision-G4.14}
\begin{equation}
     e_0(x,\gl)=((I+R)Y(\cdot,\gl))(x)=Y(x,\gl)+\int_0^x R(x,t) Y(t,\gl) dt.
     \label{G4-3.21}
\end{equation} 

\begin{lemma}\label{S1-1.4} Let $L$ be of class $(T_B)$.

{\rm 1.} Let 
$\sigma$ be the spectral function of the boundary value problem 
\myref{G3-2.1} and let
$g\in L^2_{\rm comp}(\R_+,\C^{2n})$. Put
 \begin{equation}
     G_0(\gl):=(\cf_{H,0}g)(\gl)=\int_0^\infty e_0(x,\gl)^* g(x) dx.
 \end{equation}
Then $G_0\in L^2_\sigma(\R)$ and if
\begin{equation}
    \int_\R G_0(\gl)^* d\sigma(\gl) G_0(\gl)=0
   \label{G1-1.20}    
\end{equation}
then $g=0$.

{\rm 2.} We have $\cf_{H,0}(L^2_{\rm comp}(\R_+,\C^{2n}))=
\cf_{H,Q}(L^2_{\rm comp}(\R_+,\C^{2n}))$.
\end{lemma}

\begin{proof} In view of   \myref{G4-3.21} we have
\begin{eqnarray}
     G_0(\gl)&=& \int_0^\infty \big[ Y(x,\gl)^* + \int_0^x Y(t,\gl)^*
     R(x,t)^*dt\big] g(x) dx\nonumber \\
    &=& \int_0^\infty Y(x,\gl)^*\big[g(x)+\int_x^\infty R(t,x)^* g(t)
     dt\big]dx,\label{G4-3.22}
\end{eqnarray}
hence $G_0(\gl)$ is also the $\fouHQ$--transform of the function
\begin{equation}
    \tilde g(x) :=((I+R^*)g)(x)=g(x) +\int_x^\infty R(t,x)^*g(t) dt.
    \label{G1-1.22}
\end{equation}

Since $g\in L^2_{\rm comp}(\R_+,\C^{2n})$ we also have
$\tilde g\in L^2_{\rm comp}(\R_+,\C^{2n})$.
This shows the inclusion 
$\cf_{H,0}(L^2_{\rm comp}(\R_+,\C^{2n}))\subset 
\cf_{H,Q}(L^2_{\rm comp}(\R_+,\C^{2n}))$. The converse inclusion
is proved analogously using \myref{mal3.3} instead of \myref{G4-3.21}.

In view of
the Parseval equality (Proposition \plref{S1-1.2})
we find
$$\int_\R G_0(\gl)^*d\sigma(\gl) G_0(\gl)=
   \int_0^\infty \tilde g(x)^*\tilde g(x)dx,$$
which by assumption \myref{G1-1.20} implies $\tilde g=0$. Since $g$ has compact
support \myref{G1-1.22} is a Volterra equation and thus $g=0$.
\end{proof}

\begin{prop}\label{S1-1.5} 
Let $\sigma$ be the spectral function of the boundary value problem 
\myref{G3-2.1} and let $\sigma_0=\frac{1}{2\pi} B_1^{-1} \gl$
be the corresponding spectral function for $Q=0$.
We abbreviate $\Sigma:=\sigma-\sigma_0$.

{\rm 1.} Let $L$ be of class $(T_B)$ and 
let $I+R$ be the transformation operator of the form
\myref{G4-3.21}. Furthermore, let $F$ be the $2n\times 2n$ matrix function defined
by \myref{G4-3.7}, i.e.
\begin{equation}
       F(x,t):=\left\{\begin{array}{ll}
        \DST  R(x,t)+\int_0^t R(x,s)R(t,s)^*ds,& x>t>0,\\[1em]
        \DST  R(t,x)^* +\int_0^x R(x,s) R(t,s)^* ds,& 0<x<t.
        	      \end{array}\right.
        \label{G4-3.23}    
\end{equation}
Then we have for all $f,g\in L^2_{\rm comp}(\R_+,\C^{2n})$
\begin{equation}
   \int_\R F_0(\gl)^*d\Sigma(\gl)G_0(\gl)=\int_0^\infty\int_0^\infty
       f(x)^*F(x,t)g(t) dx dt,
   \label{G1-1.25}
\end{equation}
where $F_0,G_0$ denote the $\cf_{H,0}$--transforms of $f,g$.

{\rm 2.} Again assuming $L$ to be of class $(T_B)$ we put
\begin{equation}
     \tilde e_0(x,\gl):=\int_0^x e_0(t,\gl) dt.
     \label{G2-1.26}
\end{equation}
Then the function 
\begin{equation}
    \tilde F(x,t):=\int_\R \tilde e_0(x,\gl)d\Sigma(\gl) \tilde e_0(t,\gl)^*
     \label{G2-1.27}
\end{equation} 
exists and has a continuous mixed 
second derivative which coincides with $F(x,t)$, i.e.
$\frac{\pl^2}{\pl x\pl t}\tilde F(x,t)=F(x,t).$

{\rm 3.} Conversely, given any
increasing $n\times n$ matrix function $\sigma$
put $\Sigma:=\sigma-\sigma_0$. If
the integral
\myref{G2-1.27} exists and has a continuous mixed second derivative $F_1(x,t):=
\frac{\pl^2}{\pl x\pl t}\tilde F(x,t) $ then
\myref{G1-1.25} holds for
all $f,g\in L^2_{\rm comp}(\R_+,\C^{2n})$ with $F_1$ instead of $F.$
\end{prop}

\begin{remark} \label{S10-3.6} 
We emphasize that 3. holds for arbitrary $L$ of the form \myref{G3-2.1}
not necessarily being of class $(T_B)$.

We note that the identity \myref{G1-1.25} characterizes the spectral
function of the problem \myref{G3-2.1}. More precisely, if $\varrho$
is an increasing (normalized) 
$n\times n$ matrix function such that \myref{G1-1.25}
holds with $\Sigma_\varrho:=\varrho-\sigma_0$ then $\varrho=\sigma$.

Indeed from \myref{G1-1.25} we infer
\begin{equation}\begin{array}{l}\DST
    \int_\R F_0(\gl)^* d\sigma(\gl) G_0(\gl)=
    \int_\R F_0(\gl)^* d\varrho(\gl)G_0(\gl),\\[1em]
    \DST F_0:=\cf_{H,0}f,\quad G_0:=\cf_{H,0}g,
   		\end{array}
    \label{G9-4.22}
\end{equation}
for all $f,g\in L^2_{\rm comp}(\R_+,\C^{2n})$.
By Theorem \plref{revision-S3.7} and Lemma \plref{S1-1.4}, 2. this implies that
\myref{G9-4.22} holds for all $F_0,G_0\in L^2_\varrho(\R)$, in particular
it holds for all $F_0, G_0\in C(\R,\C^n)$ with compact support.
Thus the vector measures $d\sigma, d\varrho$ and hence the right--continuous
functions $\varrho,\sigma$ coincide.
\end{remark}
\begin{proof} 1. In view of \myref{G4-3.22} $F_0$ is the $\fouHQ$--transform of

\begin{equation}
   \widetilde f(x)=f(x)+\int_x^\infty R(t,x)^* f(t) dt,
   \tag{\ref{G1-1.22}'}
\end{equation}
thus the Parseval equality \myref{G2-2.19} gives
\begin{eqnarray*}
    &&\int_\R F_0(\gl)^* d\Sigma(\gl) G_0(\gl)= \int_\R
   F_0(\gl)^* d\sigma(\gl)
    G_0(\gl)- (f,g)=(\tilde f,\tilde g)-(f,g)\\
   &=& \int_0^\infty (f(x)+\int_x^\infty R(t,x)^* f(t) dt)^*(g(x)+\int_x^\infty
     R(t,x)^* g(t) dt) dx -(f,g)\\
   &=& \int_0^\infty\int_0^\infty f(x)^* F(x,t) g(t) dx dt
\end{eqnarray*}
by a straightforward calculation.

2. For $x,t\ge 0$ and $f_0, g_0\in \C^{2n}$ we apply 1. with
$f(u):=1_{[0,x]}(u)f_0, g(v):=1_{[0,t]}(v)g_0$ and find
\begin{equation} f_0^* \int_R \tilde e_0(x,\gl)d\Sigma(\gl)
   \tilde e_0(t,\gl)^* g_0=
    f_0^* \int_0^x\int_0^y F(u,v) dudv g_0,
    \label{G2-1.28}
\end{equation}
which implies the first assertion.

3.
To prove the converse statement we note that now we have \myref{G2-1.28}
with $F_1(x,t)=\frac{\pl^2}{\pl x\pl t}\tilde F(x,t)$. This identity
implies
\myref{G1-1.25} with $F_1$ instead of $F$ for step functions 
\begin{equation}
 f=\sum_{j=1}^n f_j 1_{[a_j,b_j[}, \quad
 g=\sum_{j=1}^n g_j 1_{[c_j,d_j[},
 \quad f_j,g_j\in \C^{2n}.
\end{equation}
There is a slight subtlety since $\Sigma$ is not necessarily increasing.
However, we conclude from \myref{G1-1.25} and the Parseval equality
that for all step functions $f,g$
\begin{equation}
   (F_0,G_0)_{L^2_\sigma(\R)}=(f,g)_{L^2(\R_+,\C^{2n})}+
    \int_0^\infty\int_0^\infty f(x)^* F(x,t) g(t) dx dt.
  \label{G9-3.34}
\end{equation}
Since $\sigma$ is increasing the assertion now follows from the
denseness of the step functions in $L^2_{\rm comp}(\R_+,\C^{2n})$.
To complete the proof it remains to note that the equality 
$F(x,t) = F_1(x,t)$ is a consequence of
\myref{G1-1.25} and \myref{G9-3.34}.
\end{proof}

Combining Propositions \plref{S4-3.1} and \plref{S1-1.5}
one immediately obtains the following theorem.

\begin{theorem}\label{S1-1.6} 
Assume that the system \myref{G3-2.1} is of class $(T_B)$.
Let $\sigma$ be its spectral measure function and 
$\sigma_0(\gl)=\frac{1}{2\pi} B_1^{-1}\gl$.
Then with $F$ defined by \myref{G2-1.27} we have 
the Gelfand--Levitan equation
\begin{equation}
   F(x,t)+K(x,t) +\int_0^x K(x,s) F(s,t) ds =0,\quad t<x.
   \label{G1-1.27}
\end{equation}
\end{theorem}
\begin{remark} Note that by Proposition \plref{S1-1.5} 2. the function
$F$ is continuous also on the diagonal. In view of \myref{G4-3.7}
the continuity of $F$ at the diagonal implies $R(x,x)=R(x,x)^*$.
\end{remark}

\begin{proof} We present a second proof of the Gelfand--Levitan equation
based on the formula \myref{G2-1.27}
for $F$, which is similar to \cite{GelLev:DDESF} and
\cite[Chap. 12]{\Lev}. 

For $f,g\in L^2_{\rm comp}(\R_+,\C^{2n})$ we consider
\begin{equation}
   {\rm I}(f,g):= \int_\R \int_0^\infty dx \int_0^\infty dt f(x)^*
   Y(x,\gl)d\sigma(\gl) e_0(t,\gl)^* g(t).
   \label{G1-1.28}
\end{equation}
Substituting \myref{mal3.3} for $Y$ we find using the
Parseval equality and Lemma \plref{S1-1.5}
\begin{align}
   {\rm I}(&f,g)= (f,g)+\int_0^\infty\int_0^\infty f(x)^*F(x,t)g(t) dxdt
     \nonumber\\
       &+\underbrace{\int_\R \int_0^\infty dx \int_0^\infty dt
        f(x)^* \int_0^\infty K(x,s) e_0(s,\gl) ds d\sigma(\gl)
         e_0(t,\gl)^* g(t)}_{{\rm II}(f,g)},\\
  {\rm II}(&f,g)= \int_\R \int_0^\infty dx \int_0^\infty ds
      \int_0^\infty f(x)^* K(x,s)dx e_0(s,\gl) d\sigma(\gl)
         e_0(t,\gl)^* g(t).\nonumber
\end{align}
Writing $d\sigma=d\Sigma+d\sigma_0$ and using Lemma \plref{S1-1.5}
we find
\begin{eqnarray}
   {\rm II}(f,g)&=& \int_0^\infty \int_0^\infty f(x)^* K(x,t) g(t) dxdt
                    \nonumber \\
    && + \int_0^\infty \int_0^\infty \int_0^\infty 
       f(x)^* K(x,s) F(s,t) g(t) dx dt ds,
\end{eqnarray}
hence 
\begin{equation}\begin{split}
  {\rm I}(&f,g)=(f,g)+\\
    &+ \int_0^\infty \int_0^\infty
    f(x)^* \big[ F(x,y)+K(x,y)+\int_0^x K(x,t)F(t,y) dt\big] g(y) dx dy.
  		\end{split}
  \label{G1-1.31}
\end{equation}

Now if $\supp f\subset [b,\infty)$, $\supp g\subset [0,a]$, $a<b$, then
$(f,g)=0$ and
$$\int_0^\infty e_0(x,\gl)^* g(x) dx$$
is the $\fouHQ$--transform of 
$$g(x) +\int_x^\infty R(t,x)^* g(t) dt$$
which also has support in $[0,a]$, hence by the Parseval equality
${\rm I}(f,g)=0.$ This implies the assertion.\end{proof}

\section{The inverse problem}\label{sec4}

\subsection{The main result}
\begin{prop}\label{S1-1.7} Let $B=\diag(B_1,-B_2)$ be an arbitrary nonsingular
self--adjoint
matrix of signature $0$.
Let $\sigma(\gl)$ be a $n\times n$ matrix
function satisfying:
\begin{enumerate}
\renewcommand{\labelenumi}{{\rm \arabic{enumi}.}}
\item If $g\in L^2_{\rm comp}(\R_+,\C^{2n})$ and if
    $$\int_\R G_0(\gl)^* d\sigma(\gl) G_0(\gl)=0,$$
where $G_0$ is the $\cf_{H,0}$--transform of $g$, then $g=0$.
\item The function
\begin{equation}
 \tilde F(x,t):=\int_\R \tilde e_0(x,\gl)d\Sigma(\gl) \tilde e_0(t,\gl)^*
   \label{G2-2.30}
\end{equation}
with $\Sigma=\sigma-\sigma_0$
exists, and has a continuous mixed second derivative 
\begin{equation}
   F(x,t):=\frac{\pl^2}{\pl x\pl t}\tilde F(x,t).
   \tag{\ref{G2-2.30}'}
\end{equation}
\end{enumerate}
Then the Gelfand--Levitan equation \myref{G1-1.27} has a unique
continuous solution $K:\gO\longrightarrow {\rm M}(2n,\C)$.

Moreover, if $F(x,t)$ is continuously differentiable, then so is
$K(x,t)$.
\end{prop}
\begin{proof} Since for fixed $x$ equation \myref{G1-1.27} is a Fredholm equation
it suffices to show that the dual equation
\begin{equation}
   k(t)+\int_0^x k(s) F(t,s)^* ds=0,
   \label{G1-1.35}
\end{equation}
where $k:[0,x]\to M{\rm }(2n,\C)$ is square integrable, has only
the zero solution. Looking at the individual columns in \myref{G1-1.35}
it suffices to show that
\begin{equation}
   g(t)^*+\int_0^x g(s)^*F(t,s)^* ds=0, \quad g\in L^2([0,x],\C^{2n})
   \label{G1-1.36}
\end{equation}
implies $g=0$. Extending $g$ by $0$ to $\R_+$ we may consider
$g$ as an element of $L^2_{\rm comp}(\R_+,\C^{2n})$ and \myref{G1-1.36}
implies in view of 2. and Proposition \plref{S1-1.5}, 3.
\begin{eqnarray*}
   0&=&\|g\|^2+\int_0^\infty\int_0^\infty g(s)^* F(s,t) g(t) ds dt\\
     &=& \|g\|^2+\int_\R G_0(\gl)^* d\Sigma(\gl) G_0(\gl)
    = \int_\R G_0(\gl)^* d\sigma(\gl) G_0(\gl)
\end{eqnarray*}
and thus $g=0$ by 1.

The proof of $C^1$--smoothness of $K(x,t)$ is similar to that used
in \cite{LevSar:SLDO} and \cite{GelLev:DDESF} and is omitted.
\end{proof}


Next we prove the main result of this paper:

\begin{theorem}\label{S2-4.1} 
Let $B=\diag(B_1,-B_2)$ be an arbitrary nonsingular self--adjoint matrix
of signature $0$ as in \myref{G1-1.1}.
Let $\sigma(\gl)$ be an increasing (right--continuous, $\sigma(0)=0$)
$n\times n$ matrix function satisfying 
the conditions {\rm 1.} and {\rm 2.} of Proposition \plref{S1-1.7}.

Then there exists a unique continuous $2n\times 2n$ matrix
potential $Q$ satisfying \myref{mal3.2} such that the corresponding
system \myref{G3-2.1} is of class $(T_B)$
and such that $\sigma$ is its spectral measure function.
$Q(x)$ has $p$
continuous derivative iff $D^p_xD^p_tF(x,t)$ is continuous.

Conversely, if $\sigma$ is the spectral measure function of the
boundary value problem \myref{G3-2.1} of class $(T_B)$
then
the conditions {\rm 1.} and {\rm 2.} of Proposition \plref{S1-1.7} hold.
\end{theorem}

\begin{proof} The necessity was proved in Lemma \plref{S1-1.4} and
Proposition \plref{S1-1.5}.

To prove the sufficiency we assume that the conditions 1. and 2.
of Proposition \plref{S1-1.7} hold:

i) Starting with $\sigma(\gl)$ we define $\tilde F, F$ by \myref{G2-2.30}
and (\ref{G2-2.30}').
Then we consider the Gelfand--Levitan equation \myref{G1-1.27}
\begin{equation}
   \Phi(x,t):= F(x,t)+K(x,t) +\int_0^x K(x,s) F(s,t) ds =0,\quad x>t.
   \label{G3-4.4}
\end{equation}
By Proposition \plref{S1-1.7} this equation has a unique continuous solution
$K:\Omega\to {\rm M}(2n,\C)$.

Then $F$ also equals the right hand side of \myref{G4-3.7}:
namely, starting with $K$ we consider the operator $R$ of the
form \myref{G4-3.3} and introduce $F_1$ by \myref{G4-3.7}.
According to Proposition \plref{S4-3.1} $F_1$ and $K$ are connected
by equation \myref{G4-3.10}. Thus $F$ defined by \myref{G2-2.30} and $F_1$
defined by \myref{G4-3.7} satisfy the equation \myref{G3-4.4} and
therefore we infer from Proposition \plref{S4-3.1} that $F=F_1$.

We collect further properties of $F$: in view of \myref{G4-3.7} we have
\alpheqn
\begin{equation}
    F(x,t)=F(t,x)^*.\label{G7-4.5a}
\end{equation}

By continuity, the equation \myref{G4-3.7} also holds for $x=t$ and
consequently $R(x,x)$ is self--adjoint. Therefore, so is
$K(x,x)=-R(x,x)$.
Furthermore,
\begin{equation}
\pl_tF(x,t)B=-B\pl_xF(x,t),
  \label{G7-4.5b}
\end{equation}
where this equality holds in the distributional sense if $F$ is only
continuous. To see this let $f,g\in\cinfz{(0,\infty),\C^{2n}}$. In view
of \myref{G1-1.25} and \myref{G4-2.22} applied with $Q=0$ we calculate
\begin{align*}
\int_0^\infty\int_0^\infty f(x)^* \pl_t &F(x,t)B g(t) dx dt\\
    &= -i \int_0^\infty\int_0^\infty f(x)^* F(x,t)\frac 1i B \pl_t g(t) dx
    dt\\
 &=-i\int_\R (\cf_{H,0} f)(\gl)^* d\Sigma(\gl) \gl (\cf_{H,0}g)(\gl)\\
 &=-i\int_\R (\cf_{H,0}\frac 1i B f')(\gl)^* d\Sigma(\gl)
        (\cf_{H,0}g)(\gl)\\
 &= -  \int_0^\infty\int_0^\infty f(x)^* B\pl_x F(x,t) g(t) dx dt.
\end{align*}
Moreover, it follows from \myref{G2-2.30} and \myref{G1-1.11} that with
some matrix function $T(t)$ we have
\begin{equation}
     F(0,t)={I\choose H} T(t).
    \label{G7-4.5c}
\end{equation}
\reseteqn

We now define (cf. \myref{mal3.3})
\begin{equation}
  Y(x,\gl)=e_0(x,\gl)+\int_0^x K(x,t) e_0(t,\gl) dt
   \label{G3-4.1}
\end{equation}
and we will show that the properties
(\ref{G7-4.5a}-c) imply that $Y(x,\gl)$ satisfies the initial value problem
\begin{equation}
   B\frac 1i \frac{dY(x,\gl)}{dx} +Q(x) Y(x,\gl)=\gl Y(x,\gl),
\quad Y(0,\gl)=e_0(0,\gl)={I\choose H},
   \label{G3-4.2}
\end{equation}
where
\begin{equation}
    Q(x):=iBK(x,x)-iK(x,x)B.
   \label{G3-4.3}
\end{equation}
Note that since $K(x,x)$ is self--adjoint $Q(x)$ is self--adjoint, too.
Moreover, from \myref{G3-4.3} we also conclude that $Q(x)$ is
off--diagonal, i.e. $Q_{ii}=0$.

It follows from \myref{G7-4.5c} that
\begin{equation}
     F(x,0)BF(0,t)= T(x)^*[B_1-H^*B_2H]T(t)=0.
     \label{G7-4.5d}
\end{equation}
Plugging \myref{G7-4.5d}
into the Gelfand--Levitan equation \myref{G3-4.4}
gives
\begin{equation}
    K(x,0) B {I\choose H}=0 \quad \mbox{\rm for}\quad x\in[0,\infty).
  \label{G7-4.5e}
\end{equation}

ii) For the moment we assume in addition that $F$ is continuously
differentiable. Then by Proposition \plref{S1-1.7} $K$ also is
continuously differentiable. Differentiating \myref{G3-4.4} we obtain
\begin{eqnarray}
    B \pl_x\Phi(x,t)&= &B \pl_xF(x,t) + B \pl_x K(x,t) + BK(x,x)F(x,t)
   \nonumber\\
    &&+\int_0^x B \pl_xK(x,s)F(s,t) ds=0,
    \label{G3-4.5}\\
   \pl_t\Phi(x,t)B&=&\pl_t F(x,t) B +\pl_tK(x,t) B\nonumber\\
     && +\int_0^x K(x,s)\pl_t F(s,t)B
      ds  =0.\label{G3-4.6}
\end{eqnarray}

Integrating by parts and using \myref{G7-4.5b} and \myref{G7-4.5e}
we obtain
\begin{eqnarray}
   &&  \int_0^x K(x,s)\pl_tF(s,t) B ds=-\int_0^x K(x,s) B \pl_sF(s,t) ds\nonumber\\
   &=&\int_0^x \pl_s K(x,s) B F(s,t) ds - K(x,x) B F(x,t).
   \label{G3-4.8}
\end{eqnarray}

Adding up \myref{G3-4.5} and \myref{G3-4.6} and using \myref{G3-4.8}
and the Gelfand--Levitan equation \myref{G3-4.4}
we obtain
\begin{equation}\begin{array}{l}\DST
    B\pl_x K(x,t) +\pl_t K(x,t) B +iQ(x)K(x,t)\nonumber\\[1em]
   \DST+\int_0^x [B\pl_x K(x,s)+\pl_sK(x,s)B+iQ(x)K(x,s)]F(s,t) ds =0.
  \label{G3-4.9}
    		\end{array}
\end{equation}
Since the homogeneous integral equation corresponding to the
Gelfand--Levitan equation \myref{G3-4.4} has only the trivial
solution (see the proof of Proposition \plref{S1-1.7}) we infer
from \myref{G3-4.9} that
\begin{equation}
    B \pl_x K(x,t) +\pl_t K(x,t) B +iQ(x) K(x,t)=0.
   \label{G3-4.10}
\end{equation}
Since $K$ satisfies the relations \myref{G7-4.5e}, \myref{G3-4.3}
and \myref{G3-4.10} it follows from Theorem \plref{revision-S4.1}
that $Y(x,\gl)$ (cf. \myref{G3-4.1}) satisfies the initial value
problem \myref{G3-4.2}.

iii) We now assume that $F$ is just continuous.
Assume for the moment that for $\delta>0$ we have a continuously
differentiable matrix function $F^\gd:\R_+^2\to {\rm M}(2n,\C)$
with the properties:
\alpheqn
\begin{list}{}{\leftmargin2em}
\item$F^\gd$ converges to $F$ as $\gd\to 0$ uniformly on compact subsets
of $\R_+^2$.\label{G7-4.16a}\puteqnnum

\item$F^\gd$ satisfies (\ref{G7-4.5a}-c).\label{G7-4.16b}\puteqnnum
\end{list}
\reseteqn

We fix $x_0>0$. For $0<x\le x_0$ let $T_F$ be the integral operator
in $C([0,x],\C^{2n})$ defined by  $(T_Ff)(t):=\int_0^x f(s) F(s,t)ds$.
The proof of Proposition \plref{S1-1.7} shows that $-1\not\in\spec T_F$.
Thus for $\delta\le \delta_0(x_0)$ we have $-1\not\in\spec T_{F^\delta}$
and the Gelfand--Levitan equations
\begin{equation}
  ((I+T_{F^\delta})K_\delta(x,.))(t)=K_\delta(x,t)+\int_0^x K_\delta(x,s)
    F^\delta (s,t)ds=-F^\delta(x,t)
   \label{G7-4.17}
\end{equation}
have (for each fixed $x\in (0,x_0]$) unique solutions $K_\delta(x,t)$,
$(x,t)\in [0,x_0]^2$, which converge to $K$ as $\delta\to 0$ uniformly
on $[0,x_0]^2$. Since $F^\delta$ is $C^1$ it can be shown
(cf. the proof of Proposition \plref{S1-1.7}) that $K_\delta$ is
$C^1$, too.

Moreover, $K_\delta$ satisfies \myref{G7-4.5e} for $0<x\le x_0$ which
follows from \myref{G7-4.16b} and \myref{G7-4.17}. Now part ii) of this
proof shows that $K_\delta$ also satisfies \myref{G3-4.10}
with $Q_\delta(x):=i BK_\delta(x,x)-iK_\delta(x,x) B, x\in [0,x_0]$.
Hence, $K_\delta$ satisfies (\ref{mal3.4a}-c) (on $[0,x_0]^2$) and
therefore,
$$Y_\delta(x,\gl):=e_0(x,\gl)+\int_0^x K_\delta(x,t)e_0(t,\gl)dt,
\quad 0\le x\le x_0,$$
satisfies the initial value problem \myref{G3-4.2} with $Q_\delta$
instead of $Q$.

Since $F^\delta(x,t)^*=F^\delta(t,x)$ one concludes as in part i)
of this proof that $Q_\delta(x)^*=Q_\delta(x)$.

Since $K_\delta$ converges to $K$ as $\delta\to 0$ uniformly
on $[0,x_0]^2$, $Q_\delta$  converges to $Q$  uniformly on $[0,x_0]$.
Thus $Y(x,\gl)$ satisfies the initial value problem
 \myref{G3-4.2} on $[0,x_0]$. Since $x_0$ was arbitrary $Y(x,\gl)$ satisfies
\myref{G3-4.2} on $\R_+$.

It remains to prove the existence of the sequence $F^\delta$:

Let $F(x,t):=(F_{ij}(x,t))_{i,j=1}^r$ be the block--matrix representation
with respect to the orthogonal decomposition
$\C^{2n}=\moplus_{i=1}^r \C^{n_i}$.

It follows from \myref{G2-2.30} and (\ref{G2-2.30}') that
\begin{equation}
F_{ij}(x,t)=f_{ij}(\mu _ix-\mu _jt),   \qquad
f_{ij}(\xi )=H_ig(\xi )H_j^*,
\label{G8-4.17}
\end{equation}
with $\mu _i=\lambda ^{-1}_i,\;1\le i\le r,$ and
$H_1:=I_{n_1}=I_n$. Here the map $g:\R\to {\rm M}(n\times n,
 \C)$ is continuous and satisfies $g(\xi )^*=g(-\xi ).$
Therefore the maps $f_{ij}:\R\to M(n_i\times n_j,\C)$ are
continuous and satisfy $f_{ij}(\xi)^*=f_{ji}(-\xi).$
We note that if the measure $\Sigma (\lambda )$ is finite, that is
$\int _{\R}|d\Sigma (\lambda )|\in M(n,\C) $, then
$g(\xi )=\int_{\R} e^{i\lambda \xi }d\Sigma (\lambda )$.

We put
\begin{equation}\begin{split}
g^\delta(\xi) &:= \frac1{2\delta}\int_{\xi-\delta}^{\xi+\delta}g(s)ds,\\
f^\delta _{ij}(\xi )&:=H_ig^\delta(\xi)H_j^*
\ , \quad  F^\delta _{ij}(x,t):=f^\delta _{ij}(\mu _ix-\mu _jt),
		\end{split}
\label{G8-4.18}
\end{equation}
and $F^\delta (x,t):=(F^\delta _{ij}(x,t))^r_{i,j=1}$.

Obviously, $F^\delta$ is continuously differentiable and satisfies
\myref{G7-4.16a}.
It is clear from \myref{G8-4.17} that
\begin{equation}\begin{split}
g^\delta(\xi )^*&=\frac{1}{2\delta }\int^{\xi +\delta }_{\xi -\delta }
g^\delta(s)^*ds=\frac{1}{2\delta }\int ^{\xi +\delta }_{\xi -\delta
}g^\delta(-s)ds\\
&=
\frac{1}{2\delta }\int ^{-\xi +\delta }_{-\xi -\delta }g^\delta(s)ds=
g^\delta(-\xi ),
		\end{split}
\label{G8-4.19}
\end{equation}
and thus $f_{ij}^\delta(\xi)^*=f_{ji}^\delta(-\xi).$

In view of \myref{G8-4.18} and \myref{G8-4.19}
$F^\delta$
satisfies (\ref{G7-4.5a},b).
To prove the property \myref{G7-4.5c}
for $F^\delta$ we note that in
view of \myref{G8-4.17} and \myref{G8-4.18}
$F^\delta _{ij}(0,t)=f^\delta _{ij}(-\mu_jt)=H_ig^\delta (-\mu_jt)H^*_j$
and consequently
\begin{equation}
F^\delta (0,t)=(F^\delta _{ij}(0,t))^r_{i,j=1}=(H_ig^\delta (-\mu_jt)
H^*_j)^r_{i,j=1}=:{I\choose H} T^\delta(t),
\label{G8-4.20}
\end{equation}
where $T^\delta(t)=(g^\delta(-\mu_1t)H_1^*,g^\delta(-\mu_2 t)H_2^*,\ldots,
g^\delta(-\mu_r t)H_r^*).$

This proves that $F^\delta$ satisfies \myref{G7-4.5c}.
Summing up, we have proved that
$F^\delta$ satisfies (\ref{G7-4.16a},b).

iv)  Starting with an increasing $n\times n$ matrix function
$\sigma(\gl)$ satisfying the conditions 1. and 2. of Proposition
\plref{S1-1.7} we have constructed the boundary
value problem \myref{G3-2.1} resp. \myref{G3-4.2}. To complete
the proof it remains to show that $\sigma(\gl)$ is, in fact,
the spectral function for the problem \myref{G3-4.2}.

Let $\varrho(\gl)$
be the spectral function of the problem \myref{G3-4.2}.
Starting with 
$\Sigma_\varrho:=\varrho -\sigma_0$ we define $F_\varrho$ by (\ref{G2-2.30}').
Then by Theorem \plref{S1-1.6} $K$ satisfies the Gelfand-Levitan
equation \myref{G1-1.27} with $F_\varrho$.
On the other hand,
in view of \myref{G3-4.4} $K$ satisfies the Gelfand--Levitan
equation with $F$ instead of $F_\varrho$. 
From Proposition \plref{S4-3.1} we infer $F=F_\varrho$.
By Remark \plref{S10-3.6} this implies $\varrho=\sigma$.
\end{proof}

\begin{remark}
1. The case $n=1$ and $B_1=B_2=1$, i.e. the case of a $2\times 2$ Dirac
system, is due to M. Gasymov and B. Levitan  
\cite{GasLev:IPDS}, \cite[Chap. 12]{LevSar:SLDO}.
We note, however, that the proof in
\cite[Chap. 12]{LevSar:SLDO} is incomplete, since
the self-adjointness of $Q$ is not proved.

2. We also note that following Krein's method 
\cite{Kre:CATPOUC} L. Sakhnovich
\cite[Chap. 3, \S3]{Sak:FPOI} 
has obtained some (implicit) sufficient conditions for
a matrix measure to be the 
spectral function of a canonical system.
\end{remark}

\subsection{Some complements to the main result}

Next we will discuss several other criteria  
which imply conditions 1. or 2. of Proposition \plref{S1-1.7}.
For brevity, in the sequel we will address them just
as "condition 1./2.".
As in the proof of  Theorem \plref{revision-S3.7}
we denote by 
$\Lambda:L^2_\sigma(\R)\to L^2_\sigma(\R), (\Lambda g)(\gl):=\gl g(\gl)$
the operator of multiplication by $\lambda$.
Furthermore, we denote by 
$\mu_T(\lambda_0)$ the multiplicity of the spectrum of a self--adjoint
operator $T$ at the point $\lambda_0$.
We first note some simple facts:

\begin{remark} \label{S8-5.4}1.
If $\Sigma(\gl)$ is increasing then
condition 1. is trivially fulfilled.

For let $g\in  L^2_{\rm comp}(\R_+,\C^{2n})$ with
\begin{equation}\begin{array}{rcl}
         0&=&\DST\int_\R G_0(\gl)^* d\sigma(\gl) G_0(\gl)\\[1em]
        &=& \DST\int_\R G_0(\gl) d\sigma_0(\gl) G_0(\gl)+
            \int_\R G_0(\gl) d\Sigma(\gl) G_0(\gl),
		\end{array}
    \label{G8-5.14}
\end{equation}
where $G_0$ is the $\cf_{H,0}$--transform of $g$. Since
$\Sigma(\gl)$ is assumed to be increasing both summands
on the right hand side of \myref{G8-5.14} are nonnegative and hence
$0$. Then Proposition \plref{S1-1.2} implies $g=0$. 

2. Assume that the matrix measure $\Sigma (\lambda )$ is finite, i.e.
   $\int _{\R} |d\Sigma (\lambda )|\in {\rm M}(n,\C)$. Then 
condition 2. is obviously fulfilled.
\end{remark}

Recall that a subset $X\subset \R$ is said to have
finite density (cf. \cite{Lev:DZEF}) if
\begin{equation}
\limsup_{R\to \infty}\frac{1}{R}|\{x\in X\,|\, |x|\le R\}|<\infty.
   \label{revision-G5.26}
\end{equation}
Otherwise, $X$ is said to have infinite density.

\begin{prop}\label{revision-S5.4}
Let $B=(B_1,-B_2)$ be as in \myref{G1-1.1}.
For an increasing $n\times n$ matrix function $\sigma$
the condition
\begin{itemize}
\item[\textnormal{1'.}] 
The set $\supp_n(d\sigma):=
\{\lambda \in \R\,|\,\mu_\Lambda(\lambda )=n\}$
has infinite density
\end{itemize}
implies condition {\rm 1.}
\end{prop}

\begin{proof}
Let $\sigma (\lambda )=(\sigma _{ij}(\lambda ))^n_{i,j=1}$ and
$\varrho (\lambda ):=\tr\sigma (\lambda )=\sigma_{11}(\lambda )+...
+\sigma _{nn}(\lambda ).$
From the inequality
\[   |\sigma_{ij}(\gl)-\sigma_{ij}(\mu)|\le 
    \sqrt{\vphantom{\sigma_{jj}}\sigma_{ii}(\gl)-\sigma_{ii}(\mu)}
    \sqrt{\sigma_{jj}(\gl)-\sigma_{jj}(\mu)},\quad \mu<\gl,\]
we infer that $d\sigma_{ij} (\lambda )$ 
is absolutely continuous with respect to
$d\varrho (\lambda )$. Hence, by the Radon--Nikodym 
Theorem there exists a density matrix
\begin{equation}
\Phi (\lambda )=(\phi _{ij}(\lambda ))^n_{i,j=1}, 
   \label{revision-G5.22}
\end{equation}
such that
\[
\sigma (\lambda )=\begin{cases}\int_{(0,\gl]}\Phi (t)d\varrho (t),& \gl\ge
0,\\
     -\int_{(\lambda,0]}\Phi (t)d\varrho (t),& \gl<0.
   		  \end{cases}
\]
Obviously, $\Phi (\lambda )\ge 0$ and thus we have
\begin{equation}\begin{split}
\supp_n(d\sigma)&=\supp (\det\Phi d\varrho)\\
&=\Bigl\{\lambda \in \supp(d\varrho) \,\Bigm|\, \int_{\gl-\eps}^{\gl+\eps}
 \det \Phi (\lambda )d\varrho(\gl)>0 \text{ for all } \eps>0\Bigr\}.
		\end{split}
    \label{revision-G5.23}
\end{equation}
To check condition 1. let
$g\in L^2([0,b],{\C}^{2n})$
with
$$
0=\int^{\infty }_{-\infty } G_0(\lambda )^*d\sigma (\lambda ) 
G_0(\lambda )=\int ^{\infty }_{-\infty }
    G_0(\lambda)^*\Phi(\gl) G_0(\gl)d\varrho(\gl).
$$
Then we have $G_0(\gl)=0$ for all $\gl\in \supp_n(d\sigma)$.

On the other hand
$G_0(\lambda )$ is an entire (vector) function of strict order
one and hence (cf. \cite{Lev:DZEF})
either $G_0=0$ or the set of its zeros is of finite density.
But since $\supp_n(d\sigma)$ is assumed to have infinite density, we must
have $G=0$.
\end{proof}

\begin{remark} 
1.
Condition 1.' of the previous proposition is satisfied if $\supp_n(d\sigma)$
has at least one finite limit point.

2. Note that if $n=1$ then $\supp_n(d\sigma)=\supp (d\sigma)$ equals
the support of the Radon measure $d\sigma$.
\end{remark}
\begin{cor}\label{revision-S5.5}
Let $B$ be as before and $\sigma_0(\gl)=\frac{1}{2\pi} B_1^{-1}\gl$. If
\begin{enumerate}
\renewcommand{\labelenumi}{{\rm (\roman{enumi})}}
\item the measure $d\sigma$ is discrete,
\item $\supp_n(d\sigma)=\supp(d\varrho)$,
\item the measure $d\Sigma=d\sigma-d\sigma_0,$ is finite, i.e. 
 $\int _{\R} |d\Sigma (\lambda )|\in {\rm M}(n,\C)$. 
\end{enumerate}
Then $\sigma$ satisfies condition {\rm 1.'} of Proposition 
\plref{revision-S5.4} and hence condition {\rm 1}.
\end{cor}
\begin{proof} Since $d\Sigma(\gl)$ is finite we have
\begin{equation}
\sigma (\lambda )=\frac{1}{2\pi} B_1^{-1}\lambda +C_{\pm} + o(1), \qquad
\lambda \to \pm \infty, 
  \label{revision-G5.28}
\end{equation}
where $C_{\pm}:=\Sigma (\pm \infty )=C^*_{\pm}$. Hence
\begin{equation}
  \varrho(\lambda )=c_0\lambda + c_1 + o(1),\qquad  \lambda \to +\infty.
     \label{revision-G5.29}
\end{equation}
Since $d\sigma$ and hence $d\varrho$ is discrete, we infer from
\myref{revision-G5.29} that the set of discontinuities of $\varrho$
has infinite density. Hence, by (ii) the set $\supp_n(d\sigma)$ 
has infinite density.
\end{proof}

\begin{remark} Note that this proof only uses the asymptotic
relation \myref{revision-G5.28} which is slightly weaker than
the finiteness of the measure $d\Sigma$ (since the $o(1)$ need
not be of bounded variation).
\end{remark}

Finally, we give a criterion for the condition 2.

\begin{prop}\label{revision-S5.6}
The condition {\rm 2.} is fulfilled if the limit
\begin{equation}
\lim_{\Lambda\to\infty}\int_{|\gl|\le \Lambda} e_0(x,\gl) d\Sigma(\gl)
e_0(t,\gl)^*
  \label{revision-G5.27}
\end{equation}
exists locally uniformly in $x,t$. Then indeed $F(x,t)$ is given
by \myref{revision-G5.27}.

This is the case if the matrix function $\Sigma(\gl)$
is integrable with respect to Lebesgue measure and satisfies
$\lim\limits_{\gl\to\pm\infty} \Sigma(\gl)=0$.
\end{prop}
\begin{proof} If the limit \myref{revision-G5.27} exists locally
uniformly in $x,t$ then we have
\begin{align*}
   \tilde F(x,t)&:=  
     \int_\R \tilde e_0(x,\gl)d\Sigma(\gl) \tilde e_0(t,\gl)^*\\
        &=\int_0^x\int_0^t
   \lim_{\Lambda\to\infty}\int_{|\gl|\le \Lambda} e_0(x',\gl) d\Sigma(\gl)
   e_0(t',\gl)^*dx'dt' 
\end{align*}
and we reach the first assertion. 

If $\Sigma(\gl)$
is integrable with respect to Lebesgue measure and satisfies
$\lim\limits_{\gl\to\pm\infty} \Sigma(\gl)=0$, then
we apply integration by parts for Lebesgue--Stieltjes integrals
to obtain
\begin{align*}
   \int_{|\gl|\le \Lambda} &e_0(x,\gl) d\Sigma(\gl)
e_0(t,\gl)^* = e_0(x,\lambda)\Sigma(\lambda)e_0(t,\gl)
   \bigm|_{\gl=-\Lambda}^{\gl=\Lambda}\\
  &\quad-\int_{|\gl|\le\Lambda} (\pl_\gl e_0(x,\gl))\Sigma(\gl)e_0(t,\gl)
     + e_0(x,\gl)\Sigma(\gl)(\pl_\gl e_0(t,\gl)) d\gl.
\end{align*}
Since $e_0(x,\gl)$ is uniformly bounded and $\pl_\gl e_0(x,\gl)$ is
uniformly bounded for $|x|\le R$ for each $R$, we reach the conclusion.
\end{proof}

\begin{cor}\label{S11-6.6} The spectral function can be prescribed on
   an arbitrary finite interval. More precisely,
   there exists a boundary value problem \myref{G3-2.1} with
   continuous $Q$ satisfying \myref{mal3.2}
   and such that its spectral function
   $\sigma (\lambda )$ coincides on an arbitrary finite interval with a
   prescribed increasing $n\times n$ spectral measure.
\end{cor}
\begin{proof}
This follows from the fact that if $\Sigma(\gl)$ is constant outside
a compact interval then it satisfies condition 2. by
Remark \plref{S8-5.4} (or the previous proposition) and it
satisfies condition 1. by Proposition \plref{revision-S5.4}.
\end{proof}

\begin{conjecture}\label{S11-6.7} 
We conjecture that condition {\rm 1.} in Theorem
\plref{S2-4.1} is obsolete in general.
\end{conjecture}

\section{Some generalizations, comments, examples}\label{sec5}

\subsection{Generalization of the main result Theorem \ref{S2-4.1}}

Before we have investigated an operator $L$ of the form
\myref{G3-2.1} starting with the operator $L_0$ (with $Q=0$).
This has an obvious generalization. Namely, we may investigate
two operators $L_1:=L_{1,H}, L_2:=L_{2,H}$ and consider $L_2$
as a perturbation of $L_1$. More precisely,  let
\begin{equation}
 L_j=\frac 1i B\frac{d}{dx}+Q_j,
  \label{G6-5.1}
\end{equation}
and
\begin{equation}
  \cd(L_j)=\{ f\in\cd(L_j^*)\,|\, f_2(0)=Hf_1(0)\},\quad
   B_1=H^*B_2H.
   \label{G6-5.2}
\end{equation}
Furthermore, let $Y_j$ be the $2n\times n$ matrix solution of the
initial value problem \myref{G1-1.16}
(with $L_j$ instead of $L$). 
If both operators $L_j$ are of class $(T_B)$ then
$Y_j$ admits the representation
$Y_j(.,\lambda )=(I+K_j)e_0(.,\lambda)$, where $K_j$ is a Volterra operator
with kernel $K_j(x,t)$.
Therefore
\begin{equation}
 Y_2(x,\lambda )=
 ((I+K)Y_1(.,\lambda ))(x)=Y_1(x,\lambda )+\int ^x_0K(x,t)Y_1(t,\lambda )dt,
 \label{G6-5.3}
\end{equation}
where 
\begin{equation}
 I+K=(1+K_2)(I+K_1)^{-1}.
 \label{G8-5.4}
\end{equation}

Repeating the arguments used in the proof of Theorem \plref{revision-S4.1}
one concludes
that if $Q_1,Q_2 \in C^1(\R_+,{\rm M}(2n,\C))$ then
$K\in C^1(\Omega, {\rm M}(2n,\C))$ and, moreover, $K$ satisfies the
following Goursat problem
\alpheqn[G7-5.4]
\begin{eqnarray}
&&B\pl_xK(x,t)+\pl_tK(x,t)B+iQ_2(x)K(x,t)-iK(x,t)Q_1(t)=0, \label{G7-5.4a}\\
&&  BK(x,x)-K(x,x)B=i(Q_1(x)-Q_2(x)),\label{G7-5.4b}\\
&&K(x,0)B{I\choose H}=0.\label{G7-5.4c}
\end{eqnarray}
\reseteqn
We also note that \myref{G7-5.4a}--\myref{G7-5.4c}
may be deduced directly from \myref{G8-5.4} and \myref{G4-3.15}
for $K_1, K_2$. For example \myref{G7-5.4b} follows from
\myref{mal3.4b} and the
identity $K(x,x)=K_2(x,x)-K_1(x,x)$.

Putting $R:=(I+K)^{-1}-I$ we obtain from \myref{G6-5.3}
\begin{equation}
 Y_1(x,\lambda )=((I+R)Y_2(.,\lambda))(x)
  =Y_2(x,\lambda )+\int ^x_0R(x,t)Y_2(t,\lambda )dt.
 \label{G6-5.4}
\end{equation}

Since Proposition \plref{S4-3.1} remains valid in the case
under consideration, the following result, being a complete analog
of Proposition \plref{S1-1.5},
may be obtained in the same way as Proposition \plref{S1-1.5}.

\begin{prop}\label{S6-5.1}
Let $\sigma_j(\lambda )$ be the $n\times n$ spectral
function (cf.  Proposition \plref{S1-1.2}) of the operator $L_j, j=1,2$, and
$\Sigma :=\sigma _2-\sigma _1$.

{\rm 1.} Let $L_j$ be of class $(T_B)$ and
let $F(x,t)$ be defined by \myref{G4-3.23} with $R(x,t)$ being
the kernel of the transformation operator \myref{G6-5.4}.
Then we have for all
     $f,g\in L^2_{\rm comp} (\R_+,\C^{2n})$
\begin{equation}
 \int_\R F_1(\lambda )^*d \Sigma(\lambda )G_1(\lambda )
=\int^\infty_0\int ^\infty_0
  f(x)^*F(x,t)g(t) dxdt,
 \label{G6-5.5}
\end{equation}
where $F_1$ and $G_1$ are the $\cf_{H,Q_1}$-transforms of $f$ and
 $g$ respectively.

{\rm 2.} Again assuming $L_j$ to be of class $(T_B)$ we put
    $$
    \widetilde Y_1(x,\lambda ):= \int ^x_0 Y_1(t,\lambda )dt.
    $$
Then the function
\begin{equation}
  \widetilde F (x,t):= \int _\R \widetilde
  Y_1 (x, \lambda )d\Sigma (\lambda ) \widetilde Y_1(t,\lambda )^*
  \label{G6-5.6}
\end{equation}
exists and has a continuous mixed second derivative which coincides
with $F(x,t)$, i.e. $\frac{\pl^2}{\pl x\pl t}\widetilde F(x,t)=F(x,t).$

{\rm 3.}
Conversely, given any increasing $n\times n$ matrix
function $\sigma _2$ put $\Sigma :=\sigma_2-\sigma_1$. If the
integral \myref{G6-5.6} exists and has a continuous mixed second derivative
 $F_1(x,t):=\frac{\pl^2}{\pl x\pl t}\widetilde F(x,t)$
then \myref{G6-5.5} holds for all
$f,g\in L^2_{\rm comp}(\R_+,\C^{2n})$ with $F_1$ instead of $F$.
\end{prop}

Again, we emphasize that 3. holds for arbitrary 
$L$ of the form \myref{G3-2.1}
not necessarily being of class $(T_B)$.

Combining Propositions \plref{S6-5.1} and
\plref{S4-3.1} we arrive at the Gelfand-Levitan equation:

\begin{prop}\label{S7-5.2} Let $L_j$ be of class $(T_B)$ and
let $\sigma_j$ be the spectral function of the problem
\myref{G3-2.1} with $Q_j=Q^*_j,$  $j=1,2$,
instead of $Q$. Then with $F$ defined
by \myref{G6-5.6} we have the Gelfand-Levitan equation \myref{G1-1.27}.
\end{prop}

Now we are ready to present a generalization of
the main result (Theorem \plref{S2-4.1}).

\begin{theorem}\label{S7-5.3}
Let $\sigma _1(\lambda)$ be the spectral function
of the operator $L_1$ of the form \myref{G6-5.1}.
For an increasing
$n\times n$ matrix function $\sigma (\lambda )$
to be the spectral function of the boundary value problem
\myref{G3-2.1} with (unique)
continuous $2n\times 2n$ matrix potential $Q$ satisfying
\myref{mal3.2} it is
sufficient that the following conditions hold:
\begin{enumerate}
\renewcommand{\labelenumi}{{\rm \arabic{enumi}.}}
\item If $g\in L^2_{\rm comp}(\R_+,\C^{2n})$ and if
    $$\int_\R G_1(\gl)^* d\sigma(\gl) G_1(\gl)=0,$$
where $G_1$ is the $\cf_{H,Q_1}$--transform of $g$, then $g=0$.
\item The function
\begin{equation}
        \widetilde F(x,t):=\int _\R \widetilde Y_1(x,\lambda )
        d\Sigma (\lambda )\widetilde Y_1(t,\lambda)^*,\qquad
         \widetilde Y_1(x,\lambda ):=
         \int ^x_0Y_1(t,\lambda )dt,
        \label{G7-5.6}
\end{equation}
with $\Sigma=\sigma-\sigma_1$
exists and has a continuous mixed second derivative
\begin{equation}
   F(x,t):=\frac{\pl^2}{\pl x\pl t}\tilde F(x,t).
   \tag{\ref{G7-5.6}'}
\end{equation}
\end{enumerate}
Moreover $Q$ has $m$ continuous derivatives if and only if 
$D^m_xD_t^mF(x,t)$ exists and is continuous.

Again, if we content ourselves to operators of class $(T_B)$ 
then the conditions {\rm 1.} and {\rm 2.} are
also necessary.
\end{theorem}
\begin{proof}[Sketch of Proof]
The necessity is proved in just the same way as Lemma \plref{S1-1.4}
and Proposition \plref{S1-1.5}.

\noindent \emph{Sufficiency:}\quad 
Starting with $\sigma(\lambda)$ we define
$\widetilde F, F$ by \myref{G6-5.6} with $\Sigma(\lambda):=\sigma(\lambda)-\sigma_1(\lambda).$
Then we consider the Gelfand-Levitan equation
\begin{equation}
F(x,t)+K(x,t)+\int ^x_0K(x,s)F(s,t)ds=0,\quad  t<x.
\label{G9-5.8}
\end{equation}
with $F$ defined by \myref{G6-5.6}. 
Following the proof of Proposition \plref{S1-1.7}
one concludes
that \myref{G9-5.8} has a continuous solution $K:\Omega \to {\rm M}(2n,\C)$.
Next we define $Y(x,\lambda )$ setting $Y(.,\lambda )=(I+K)Y_1(.,\lambda )$
and show that $Y(x,\lambda )$ satisfies the initial value problem \myref{G3-4.2} with
\begin{equation}
Q(x)=Q_1(x)+iBK(x,x)-iK(x,x)B.
\label{G9-5.9}
\end{equation}
Since $Q_1$ satisfies \myref{mal3.2} we infer from \myref{G9-5.9}
that $Q$ also satisfies \myref{mal3.2}.
Moreover the self-adjointness of $Q$ may be proved as in the proof of Theorem 
\plref{S1-1.7}.

Furthermore,  we note that if $F$ is continuously differentiable it
satisfies the equality
\begin{equation}
BD_xF(x,t)+D_tF(x,t)B=-iQ_1(x)F(x,t)+iF(x,t)Q_1(t)
   \label{G9-5.10}
\end{equation}
and according to Proposition \plref{S1-1.7} $K$ is continuously 
differentiable, too. If $F$ is just continuous then \myref{G9-5.10}
still holds in the distributional sense. This is shown similar
to \myref{G7-4.5b}.

Since $Y_1(0,\lambda )={I\choose H}$ we may argue exactly as in 
\myref{G7-4.5c}, \myref{G7-4.5d}, \myref{G7-4.5e}
to obtain
\begin{equation}
K(x,0)B {I\choose H}=0,  \quad \mbox{\rm for}\quad x\in[0,\infty).
\label{G9-5.11}
\end{equation}

In view of \myref{G9-5.8}-\myref{G9-5.11} the relation
\myref{G7-5.4a} for $K$ is proved along the same lines
as part ii) of the proof
of Theorem \plref{S2-4.1}.

Thus $K$ satisfies the
initial value problem \myref{G7-5.4a}--\myref{G7-5.4c}. 
Therefore $Y(x,\lambda )$ satisfies
the initial value problem \myref{G3-4.2} with $Q$ defined by \myref{G9-5.9}.

If now $F$ is just continuous then one proceeds as in 
part iii) of the proof of Theorem \plref{S2-4.1}.

That $\sigma$ is indeed the spectral function of the problem
\myref{G3-4.2} with $Q$ from \myref{G9-5.9} is shown as part iv)
of Theorem \plref{S2-4.1}. Instead of 
\myref{G2-2.30},
Theorem \plref{S1-1.6},
\myref{G3-4.4},
and Proposition \plref{S1-1.5}
one uses
 \myref{G6-5.6},
 Proposition \plref{S7-5.2},
 \myref{G9-5.8},
and Proposition \plref{S6-5.1}.
\end{proof}

\subsection{The degenerate Gelfand--Levitan equation}

We discuss solutions of the Gelfand--Levitan equation in the special case where
$\Sigma(\gl)$ is a step function:

We consider the situation of Theorem \plref{S7-5.3} and fix an operator $L_1$
of the form \myref{G6-5.1} with spectral function $\sigma_1(\gl)$.

Let $A\in {\rm M}(n,\C)$ be a hermitian nonnegative matrix and
\begin{equation}
    \Sigma(\gl):= A \;1_{[a,\infty)}(\gl)
   \label{G10-5.15}
\end{equation}
an increasing step function with one jump of ``height'' $A$.

We show that
\begin{equation}
    \sigma:= \sigma_1+\Sigma
\end{equation}
is the spectral function of the boundary value problem \myref{G3-2.1}
for some (unique) continuous self--adjoint $2n\times 2n$--matrix potential
$Q$ satisfying \myref{mal3.2}.

Since jumps of the spectral function correspond to eigenvalues
this shows in particular that for a given potential
$Q_1$ and given real number $a$ 
there is a potential $Q$ such that
\begin{equation}
      {\rm spec}(L_1+Q-Q_1)={\rm spec}(L_1)\cup \{a\}.
\end{equation}

For the proof we have to verify the conditions 1. and 2. of Theorem 5.3.
By Remark \plref{S8-5.4} condition 1. is fulfilled since $A$ is
nonnegative. To verify 2. we calculate
\begin{eqnarray*}
       \tilde F(x,t)&=& \DST \int_\R \widetilde Y_1(x,\gl)d\Sigma(\gl) \widetilde
                        Y_1(t,\gl)^*\\
     &=&  \widetilde Y_1(x,a) A \widetilde Y_1(t,a)^*.
\end{eqnarray*}
Obviously, this has a continuous mixed second derivative, namely
\begin{equation}
     F(x,t):= \frac{\pl^2}{\pl x \pl t} \tilde F(x,t) = Y_1(x,a)AY_1(t,a)^*.
\end{equation}

In this case we can solve the Gelfand--Levitan equation explicitly.
First we introduce for $x>0$
\begin{equation}
      T(x):= \int_0^x Y_1(s,a)^*Y_1(s,a) ds .
\end{equation}
From 
$$ Y_1(0,a)^*Y_1(0,a)=I+H^*H\ge I$$
we infer that $T(x)>0$ is positive definite for $x\ge 0$.

We put for $t\le x$
\begin{equation}\begin{array}{rcl}
    K(x,t)&:=& -Y_1(x,a) A Y_1(t,a)^*+\\
      &&\quad +Y_1(x,a)A
                T(x) A^{1/2}(I+A^{1/2}T(x)A^{1/2})^{-1}A^{1/2}Y_1(t,a)^*\\[1em]
    &=& -Y_1(x,a) A^{1/2}(I+A^{1/2}T(x)A^{1/2})^{-1}A^{1/2} Y_1(t,a)^*.
		\end{array}
   \label{G9-5.20}
\end{equation}
Note that $(I+A^{1/2}T(x)A^{1/2})\ge I$ is positive definite, 
thus invertible. We abbreviate
$S(x):=A^{1/2}(I+A^{1/2}T(x)A^{1/2})^{-1}A^{1/2}.$ 
If $A$ is positive definite then we simply have 
$S(x)=(A^{-1}+T(x))^{-1}$.

One immediately checks that $K(x,t)$ solves the Gelfand--Levitan equation
\myref{G9-5.8} corresponding to $F$ and consequently determines $Q$
by means of of \myref{G9-5.9}.

Summarizing the previous considerations we arrive at the
   following proposition.

\begin{prop}\label{S10-5.5} Let $L_1$ be an operator of the form \myref{G6-5.1},
\myref{G6-5.2} with the spectral function $\sigma _1(\lambda )$ and let
   $\Sigma (\lambda )$ be of the form \myref{G10-5.15}.
   Then $\sigma =\sigma _1+\Sigma $ is the spectral function of
   the boundary value problem \myref{G3-2.1}
   with $2n\times 2n$ matrix potential
   $$
   Q(x)=   
Q_1(x)+i\{Y_1(x,a)S(x)Y_1^*(x,a)B-BY_1(x,a)S(x)Y_1^*(x,a)\}.
   $$
\end{prop}

\begin{cor}\label{S10-5.6}
Under the assumptions of the previous Proposition \plref{S10-5.5} let
   $Q_1=0$ (i.e. $L_1=-iB\frac{d}{dx}$). 
Then  the $2n\times 2n$ matrix potential corresponding to the spectral
function  $\sigma (\lambda )=\frac{1}{2\pi\gl_1 }\lambda I_n+\Sigma (\lambda )$
   with one jump of "height" $A$ is given by
   $$
   Q(x)=ie^{iaB^{-1}x}\{\tilde S(x)B-B\tilde S(x)\}
   e^{-iaB^{-1}x},
   $$
where 
$\tilde S(x):={I\choose H} A^{1/2}(I+xA^{1/2}(I+H^*H)A^{1/2})^{-1}A^{1/2}(I,H^*)$.
\end{cor}

If $\Sigma$ is a general increasing step function then the Gelfand--Levitan
equation is still solvable. However, we do not have such an explicit formula:

\begin{prop}\label{S8-5.5} Let $L_1$ be an operator of the form \myref{G6-5.1},
\myref{G6-5.2}
with spectral function $\sigma_1(\gl)$. Furthermore, let
$-\infty<a_1<\ldots< a_r<\infty$ be real numbers and
$A_j\in{\rm M}(n,\C)$ nonnegative matrices. 

Then for the increasing step function
$$\Sigma(\gl):= \sum_{j=1}^r A_j 1_{[a_j,\infty)}(\gl)$$
there exists a unique continuous matrix potential 
$Q$ satisfying \myref{mal3.2}
such that $\sigma_1+\Sigma$ is the spectral function of $L_1+Q-Q_1$.

Namely, $Q$ is uniquely determined by \myref{G9-5.9} with $K(x,t)$
being the solution of the Gelfand--Levitan equation \myref{G9-5.8}.
\end{prop}
\begin{proof} This follows by induction from Theorem \plref{S7-5.3} and the preceding
discussion.

The conditions 1. and 2. of Theorem \plref{S7-5.3}
can also immediately be checked directly:
namely, condition 1. is fulfilled in view of Remark \plref{S8-5.4}
since $\Sigma$ is increasing. Condition 2. follows immediately from
\[  \tilde F(x,t)=\sum_{j=1}^r  \widetilde Y_1(x,a_j) A_j \widetilde
Y_1(t,a_j)^* \quad\mbox{\rm and}\quad
F(x,t)=\sum_{j=1}^r  Y_1(x,a_j) A_j Y_1(t,a_j)^*.\]
\end{proof}

\subsection
{On unitary invariants of $2n\times 2n$ systems}
It is well--known that a selfadjoint operator $A$ in a Hilbert space is
uniquely determined (up to unitary equivalence) by the spectral type $[E_A]$
and the multiplicity function $N_{E_A}$. In this section
we will show that there exist potentials $Q$ such that
the corresponding operator $L_H$ has constant multiplicity one and
$[E]$ is of pure type (absolute continuous, singular continuous,
pure point).

\begin{dfn}
An increasing function $\mu:\R\to\R$ on the real line will be called 
$p$--admissible
if there exists a strictly increasing sequence of
real numbers, $(x_\nu)_{\nu\in\Z}$, such that
\begin{enumerate}
\item $x_0=0$,
\item the sequence $(x_{\nu+1}-x_\nu)$ is square summable,
\item $\lim\limits_{\nu\to\pm\infty} x_\nu=\pm\infty$,
\item $\mu(a_{\nu j})<\mu(a_{\nu,j+1}),$ where 
$a_{\nu j}:= x_\nu+j2^{-n-p}(x_{\nu+1}-x_\nu),$\\
$j=0,\ldots,2^{n+p}-1.$
\end{enumerate}
\end{dfn}

In particular, a strictly increasing function $\mu$ is $p$--admissible
for any $p$.

We will show that for a $p$--admissible increasing function
$\mu $ there exists an operator $L_H$ of the form
\myref{G3-2.1} such that its spectral measure $E:=E_{L_H}$ satisfies
$$
[E]=[d\mu ], \qquad  N_E(x)=p \text{ for } \mu\text{- a.e. } \ x\in {\R}.
$$

In particular there exist $2n\times 2n$ systems such that each point in
the spectrum has multiplicity one.
To prove this result we will use the criteria from the end of Section 5.

\begin{prop}\label{revision-S6.7} Let $B=(B_1,-B_2)\in {\rm M}(2n,{\C})$ 
be a matrix as in \myref{G1-1.1}
and  $H\in {\rm M}(n,\C)$ as in \myref{G1-1.8}.
Let $\mu$ be a $p$--admissible increasing function on the real line,
$1\le p\le n$.
Then there exists a continuous
potential $Q$ satisfying \myref{mal3.2} and 
such that the corresponding operator $L_H$
is unitary equivalent to the operator 
$\Lambda _p=\oplus _1^p \Lambda _1$, where
$$
\Lambda _1: L^2_{\mu }({\R})\to L^2_{\mu }({\R}),\quad
\Lambda _1f(\lambda )=\lambda f(\lambda ).
$$
\end{prop}
\begin{proof}
Let $(\psi_j)$ be the Rademacher functions
\cite[Sec. I.3]{Zyg:TS}, i.e. $\psi_1:\R\to \R$ is a function of period one,
such that
\begin{align}
      \psi_1(x)&=\begin{cases} \hphantom{-}1,& 0< x\le\frac 12,\\
                              -1,&            \frac 12< x\le 1,
                \end{cases}\\
\intertext{and}
      \psi_j(x)&=\psi_1(2^jx).
\end{align}
(In \cite[Sec. I.3]{Zyg:TS} one puts $\psi_1(1/2)=0$).
$\psi_j$ takes values $\pm 1$. The set $(\psi_j)$ is orthonormal
in $L^2[0,1]$ and
\begin{equation}
    \int_0^1 \psi_j(x) dx =0.
   \label{revision-6.22}
\end{equation}
We put
\begin{align}
\psi_0(\gl)&:=
\left(
\begin{array}{cccc}
\psi_1(\gl)&\psi_2(\gl)&\ldots&\psi_n(\gl)\\
\psi_2(\gl)&\psi_3(\gl)&\ldots&\psi_{n+1}(\gl)\\
\hdotsfor4\\
\psi_p(\gl)&\psi_{p+1}(\gl)&\ldots&\psi_{n+p-1}(\gl)
\end{array}
\right),\\
  \psi(\gl)&:=\psi_0(\frac{\gl-x_\nu}{x_{\nu+1}-x_\nu}), \quad\text{for }
        x_\nu\le \gl< x_{\nu+1},\\
    \Phi(\gl)&:=p^{-1}\frac{1}{2\pi} B_1^{-1/2}\psi(\gl)^*\psi(\gl)B_1^{-1/2}.
\end{align}
$\Phi(\gl)$ is a symmetric nonnegative
matrix of rank $p$ for each $\gl\in\R$.

Furthermore, let 
\begin{equation}
  f(x):=\frac{a_{\nu,j+1}-a_{\nu j}}{\mu(a_{\nu,j+1})-\mu(a_{\nu j})},
   \quad a_{\nu j}<x\le a_{\nu,j+1},
\end{equation}
and put 
\begin{equation}
   \varrho(\gl):=\begin{cases} \int_{(0,\gl]} f(t)d\mu(t),&\gl>0,\\
                              -\int_{(\gl,0]} f(t)d\mu(t),&\gl\le 0.
   		 \end{cases}
\end{equation}
Obviously the measures $d\varrho$  and $d\mu$ 
are mutually equivalent and 
\begin{equation}
   \varrho(a_{\nu j})=a_{\nu j}.
\end{equation}

Finally we put
\begin{equation}
\sigma(\gl):=\begin{cases} \int_{(0,\gl]} \Phi(t)d\varrho(t),&\gl>0,\\
                              -\int_{(\gl,0]} \Phi(t)d\varrho(t),&\gl\le 0.
   		 \end{cases}
\end{equation}

Note that in view of \myref{revision-6.22} and the orthonormality
of the Rademacher functions we have
$\sigma(x_\nu)=\frac{1}{2\pi} B_1^{-1}x_\nu$.
Again, by the orthonormality of the Rademacher functions and
the fact that the entries of $\psi$ are constant on the intervals
$(a_{\nu j},a_{\nu,j+1}]$,
we have for $x_\nu< \gl\le x_{\nu+1}$
\begin{equation}
   \|\sigma(\gl)- \frac{1}{2\pi}B_1^{-1}x_\nu \|\le
      \int_{(x_\nu,\gl]} \|\Phi(t)\| d\varrho (t)\le C (x_{\nu+1}-x_\nu),
\end{equation}
thus
\begin{equation}\begin{split}
      \int_{x_\nu}^{x_{\nu+1}}
    \|\sigma(\gl)- \frac{1}{2\pi}B_1^{-1}\gl \|d\gl
       &\le C' \int_{x_\nu}^{x_{\nu+1}}  (x_{\nu+1}-x_\nu)d\gl\\
       &=C' (x_{\nu+1}-x_\nu)^2.
                \end{split}
\end{equation}
Since $(x_{\nu+1}-x_\nu)_\nu$ is square summable, we infer that
the function
$$ \Sigma(\gl):=\sigma(\gl)-\sigma_0(\gl)$$
is integrable with respect to Lebesgue measure and $\lim\limits_{\gl\to
\pm\infty}\Sigma(\gl)=0.$ In view of Proposition
\plref{revision-S5.6} it satisfies condition 2. 

To show
that it satisfies condition 1. let
$g\in L^2([0,b],\C^{2n})$ with
\[0=\int_{-\infty}^\infty G_0(\gl)^*d\sigma(\gl)G_0(\gl)=
    \int_{-\infty}^\infty G_0(\gl)^* \Phi(\gl) G_0(\gl)d\varrho(\gl).\]
In view of 4. of the definition of admissibility we infer that
in each interval $(x_\nu,x_{\nu+1}]$ there exist points
$\gl_{\nu j}, j=0,..,2^{n+p}-1$ ,
such that the vector $B_1^{-1/2}G_0(\gl_{\nu j})$
lies in the null space of the matrix $\psi_0(j2^{-n-p})$.
Since $(x_{\nu+1}-x_\nu)$ is square summable,
each of the sequences $(\gl_{\nu j})_\nu$ is a sequence of infinite
density. 
Noting that
$G_0(\lambda )$ is an entire (vector) function of strict order
one and of finite type we infer as in the proof
of Proposition \plref{revision-S5.4} that $B_1^{-1/2}G_0(\gl)$ lies
in the null space of the matrix $\psi_0(j2^{-n-p}), j=0,...,2^{n+p}-1$,
for each $\gl$. It is easy to check that the intersection of
these null spaces is $0$, hence $G_0(\gl)=0$.

By Theorem \plref{S2-4.1} there exists $Q$ satisfying
\myref{mal3.2} such that $\sigma$ is the spectral measure
function of the corresponding operator $L_H$. This proves the theorem.

\end{proof}

\begin{cor} For $1\le p\le n$ 
there exist continuous potentials $Q:\R_+\longrightarrow
{\rm M}(n,\C)$ satisfying \myref{mal3.2} such that 
the corresponding operator $L_H$ has 
\begin{enumerate}
\renewcommand{\labelenumi}{{\rm \roman{enumi})}}
\item absolute continuous spectrum of multiplicity $p$,
\item singular continuous spectrum of multiplicity $p$,
\item pure point spectrum of multiplicity $p$.
\end{enumerate}
\end{cor}
\begin{proof}
We only have to note that there exist admissible increasing
functions $\mu_{\rm ac}, \mu_{\rm sc}, \mu_{\rm pp}$ such
that the measures $d\mu_{\rm ac}, d\mu_{\rm sc}, d\mu_{\rm pp}$ 
are absolute continuous, singular continuous, discrete, respectively.
Such measures obviously exist.
\end{proof}

\end{document}